\documentclass[12pt]{amsart}
\usepackage{amssymb}  
\usepackage{latexsym} 
\usepackage[all]{xy}

\newcommand{\charic}{{\operatorname{char}\,}} 
\newcommand{\GL}{\operatorname{GL}}
\newcommand{\Aut}{\operatorname{Aut}}
\newcommand{\Jac}{\operatorname{Jac}}
\newcommand{\adj}{\operatorname{adj}}
\newcommand{\eps}{{\varepsilon}}
\newcommand{\GGamma}{{\widetilde{\Gamma}}}
\newcommand{\isom}{{\, \cong \,}}

\newcommand{\PP}{{\mathbb P}}
\newcommand{\Q}{{\mathbb Q}}
\newcommand{\Gm}{{\mathbb G}_m}
\newcommand{\Z}{{\mathbb Z}}
\newcommand{\Pol}{{\operatorname{Pol}}}
\newcommand{\Hom}{{\operatorname{Hom}}}
\newcommand{\Tr}{{\operatorname{Tr}}}

\newcommand{\x}{{\bf x}}
\newcommand{\cdual}{{\mathfrak c}}
\newcommand{\DD}{{\mathfrak D}}
\newcommand{\cc}{{\mathbf c}}
\newcommand{\D}{{\mathbf D}}
\newcommand{\Kbar}{{\overline{K}}}
\newcommand{\Gal}{{\operatorname{Gal}}}
\newcommand{\SL}{{\operatorname{SL}}}
\newcommand{\PSL}{{\operatorname{PSL}}}
\newcommand{\Diag}{{\operatorname{Diag}}}
\newcommand{\PGL}{{\operatorname{PGL}}}
\newcommand{\ra}{{\longrightarrow}}
\newcommand{\im}{\operatorname{im}}
\newcommand{\G}{{\mathcal G}}
\newcommand{\Sha}{\mbox{\wncyr Sh}}

\newfont{\wncyr}{wncyr10 at 12pt}
\newfont{\wncyrten}{wncyr10 at 10pt}

\newenvironment{Proof}{\par\noindent{\sc Proof:}}%
                      {\hspace*{\fill}\nobreak$\Box$\par\medskip}
\newenvironment{ProofOf}[1]{\par\noindent{\sc Proof of #1:}}%
                      {\hspace*{\fill}\nobreak$\Box$\par\medskip}

\newtheorem{Proposition}{Proposition}[section]
\newtheorem{Theorem}[Proposition]{Theorem}
\newtheorem{Lemma}[Proposition]{Lemma}
\newtheorem{Corollary}[Proposition]{Corollary}

\theoremstyle{definition}

\newtheorem{Definition}[Proposition]{Definition}
\newtheorem{Remark}[Proposition]{Remark}

\numberwithin{equation}{section}

\newcounter{nootje}
\begin{document}

\date{25th November 2010}
\title[The Hessian of a genus one curve]
{The Hessian of a genus one curve} 
\author{Tom Fisher}
\address{University of Cambridge,
        DPMMS, Centre for Mathematical Sciences,
         Wilberforce Road, Cambridge CB3 0WB, UK}
 \email{T.A.Fisher@dpmms.cam.ac.uk}

\begin{abstract}
We continue our development of the invariant theory
of genus one curves with the aim of computing
certain twists of the universal family of elliptic curves 
parametrised by the modular curve $X(n)$ for $n=2,3,4,5$. 
Our construction makes use of a covariant we call
the Hessian, generalising the classical Hessian
that exists in degrees $2$ and $3$. In particular we 
give explicit formulae and 
algorithms for computing the Hessian in degrees $4$ and $5$.
This leads to a practical algorithm for computing equations for 
visible elements of order $n$ in the Tate-Shafarevich 
group of an elliptic curve. Taking Jacobians we 
also recover the formulae of Rubin and Silverberg for 
families of $n$-congruent elliptic curves. 
\end{abstract}

\maketitle

\renewcommand{\baselinestretch}{1.1}
\renewcommand{\arraystretch}{1.3}

\renewcommand{\theenumi}{\roman{enumi}}


\vspace{-3ex}

\section{Introduction}
\label{intro}
In our earlier paper \cite{g1inv} we developed the 
invariant theory of genus one curves of 
degrees $n=2,3,4,5$, our main original
contribution being in the case $n=5$. In this paper we
study the covariants of a genus one curve, and find that
the classical Hessian in degrees $2$ and $3$
has a natural generalisation to degrees $4$ and $5$. 
The arithmetic significance of the Hessian is that it allows us to
compute certain twists of the universal family of elliptic curves
parametrised by the modular curve $X(n)$. The existence of the
Hessian is most easily shown by a double application of 
the evectant construction, described by Salmon \cite[Art. 233]{Salmon} 
in the case $n=3$. For this reason
we study the contravariants in parallel with the covariants. 
In the case $n=5$ taking evectants is not practical since the 
invariants are too large to write down as explicit polynomials.
Nonetheless we have found a practical algorithm
for evaluating the Hessian in this case.

In Sections~\ref{sec:discinv} and \ref{sec:hesse} 
we relate the invariants of a genus one curve
to the invariants described by Klein in his {\em Lectures on the
icosahedron} \cite{Kicosa}. We extend these methods 
in Sections~\ref{sec:disccov} and \ref{sec:covcontra} 
to show that the covariants and contravariants
each form a free module of rank $2$ over the ring of invariants.
The invariants, covariants and contravariants are related by
identities recorded in Sections~\ref{sec:hpolys} and \ref{sec:hpolys2}. 
We call the polynomials arising in this context the {\em Hesse polynomials}.
In Section~\ref{sec:formulae} we give formulae for the Hessian and 
the contravariants in the cases $n=2,3,4$. Our algorithm for evaluating
the Hessian in the case $n=5$ is described in Section~\ref{sec:evaluation}.

We then turn to arithmetic applications. In Section~\ref{Theta groups}
we discuss the relationship between a genus one normal curve 
$C \to \PP^{n-1}$ and the $n$ by $n$ matrices that describe the action of
the $n$-torsion of its Jacobian. This is useful both for the later
sections of this paper, and for the Hesse pencil method of $n$-descent, 
as described in \cite[Section 5.1]{paperI}. 

In Section~\ref{famec} we show that the Hesse polynomials
define the family of elliptic curves directly $n$-congruent to a 
given elliptic curve for $n=2,3,4,5$. These formulae were previously 
obtained by Rubin and Silverberg 
\cite{RubinSilverberg}, \cite{RubinSilverberg2}, 
\cite{Silverberg} by a different method. Let us note however that
our formulae in the case $n=3$ were already
known to Salmon \cite[Art. 230]{Salmon}. 
One advantage of our method is that it generalises immediately to
the family of elliptic curves reverse $n$-congruent to a given
elliptic curve. 

In Sections~\ref{sec:vis} and~\ref{sec:examples} we use 
the Hessian to compute equations for
elements in the Tate-Shafarevich group of an elliptic curve
that are visible in the sense of Mazur \cite{CM}, 
\cite{Mazurcubics}. In the terminology of \cite{Mazurcubics} 
we have developed the invariant theory necessary to compute both
first and second twists.

We have contributed all formulae and algorithms in this paper
to the computer algebra system {\sf MAGMA} \cite[Version 2.13]{magma}.

\section{Background and overview}
\label{sec:inv}

We work over a perfect field $K$ of characteristic not dividing $6n$. 
We write $\Kbar$ for the algebraic closure and identify all $K$-schemes 
with their sets of $\Kbar$-points.
We recall some of our notation and results from~\cite{g1inv}.


\begin{Definition}
\label{def:1.1}
A {\em genus one model} of degree $n=2,3,4,5$ is \\
(i) if $n=2$ a binary quartic \\
(ii) if $n=3$ a ternary cubic \\
(iii) if $n=4$ a pair of quadrics in $4$ variables \\
(iv) if $n=5$ a $5 \times 5$ alternating matrix of 
linear forms in $5$ variables.
\end{Definition}
We write $X_n$ for the space
of genus one models of degree $n$. It is an 
affine space of dimension $N=10n/(6-n)$.  
We give the co-ordinate ring $K[X_n]$ its usual grading
by degree. A model $\phi \in X_n$ defines a subvariety
$C_{\phi}$ of $\PP(1,1,2)$ or $\PP^{n-1}$ according as 
$n=2$ or $n=3,4,5$. In the case $n=5$ the equations
are the $4 \times 4$ Pfaffians of $\phi$.
A model $\phi$ is non-singular if $C_{\phi}$ is 
a smooth curve of genus one.

In \cite{g1inv}
we defined a linear algebraic group $\G_n$ acting on $X_n$ 
and said that models $\phi, \phi' \in X_n$ are equivalent if
they belong to the same $\G_n$-orbit. We also defined a rational
character on $\G_n$ by
$$ \begin{array}{l@{\qquad}lrcl}
n = 2 & \det: \Gm \times \GL_2 \to \Gm; & [\mu,B] & \mapsto & \mu \det B\\ 
n = 3 & \det: \Gm \times \GL_3 \to \Gm; &[\mu,B] & \mapsto & \mu \det B \\ 
n = 4 & \det: \GL_2 \times \GL_4 \to \Gm; &[A,B] & \mapsto & \det A \det B \\
n = 5 & \det: \GL_5 \times \GL_5 \to \Gm; &[A,B] & \mapsto & 
(\det A)^2 \det B. 
\end{array} $$
Notice that the definitions of $\G_2$ and $X_2$ are slightly
different from those in \cite[Section~3.2]{g1inv}, since we will
not be working over fields of characteristic $2$.
We write $G_n$ for the commutator subgroup of $\G_n$.
Thus $G_2 = \SL_2$, $G_3 = \SL_3$, $G_4 = \SL_2 \times \SL_4$
and $G_5 = \SL_5 \times \SL_5$.
 

\begin{Definition}
The {\em ring of invariants} is
$$ K[X_n]^{G_n} = 
  \{ F \in K[X_n] : F \circ g = F \text{ for all } g \in G_n \}. $$
An invariant $F$ has {\em weight} $k$ if $F \circ g = (\det g)^k F$
for all $g \in \G_n$.
\end{Definition}

\begin{Lemma}
\label{degwt}
Every homogeneous invariant of degree $d$ has 
weight $k$ where $d=kn/(6-n)$.
\end{Lemma}
\begin{Proof}
This is \cite[Lemma 4.3]{g1inv}.
Some care is needed in the 
case $n=2$ since we have changed the definitions of $\G_2$ and $X_2$.
\end{Proof}

\begin{Theorem}
\label{invthm}
Let $n=2,3,4,5$. There are invariants $c_4$, $c_6$ and $\Delta$ 
of weights $4$, $6$ and $12$,
related by $c_4^3-c_6^2 = 1728 \Delta$, 
such that 
\begin{enumerate}
\item The ring of invariants $K[X_n]^{G_n}$ is 
generated by $c_4$ and $c_6$. 
\item A model $\phi \in X_n$ is non-singular 
if and only if $\Delta(\phi) \not= 0$. 
\item If $\phi \in X_n(K)$ is non-singular 
then $C_{\phi}$ is a smooth curve of genus one
defined over $K$ with Jacobian
$y^2 = x^3 - 27 c_4(\phi) x - 54 c_6(\phi)$.
\end{enumerate}
\end{Theorem}
\begin{Proof}
This is \cite[Theorem 4.4]{g1inv}. 
The cases $n=2,3,4$ are classical: 
see for example~\cite{Mc+}, \cite{Weil}.
\end{Proof}

For $g \in \G_n$ we write $g^T$ for the element obtained
by transposing the constituent matrices. We also write
$g^{-T}$ for $(g^{T})^{-1}$.

\begin{Definition} 
\label{def:cov}
A polynomial map $F : X_n \to X_n$ defined over $K$ is 
\begin{enumerate}
\item a {\em covariant} if $F \circ g = g \circ F$ for all $g \in G_n$, 
\item a {\em contravariant} if $F \circ g = g^{-T} \circ F$ 
for all $g \in G_n$.  
\end{enumerate}
A covariant or contravariant $F$ 
is homogeneous of degree~$d$ 
if $F(\lambda \phi) = \lambda^d F(\phi)$ for all $\lambda \in \Kbar$ 
and $\phi \in X_n$. 
It has {\em weight} $k$ if 
$F \circ g = (\det g)^k g \circ F$, respectively
$F \circ g = (\det g)^k g^{-T} \circ F$ for all $g \in \G_n$.
\end{Definition}

It is clear that the covariants and contravariants each
form a module over the ring of invariants $K[X_n]^{G_n}=K[c_4,c_6]$.

\begin{Lemma}
\label{prop:covwts}
Every homogeneous covariant, respectively contravariant, of degree $d$
has weight $k$ where $d=1+kn/(6-n)$, respectively $d=-1+kn/(6-n)$.
\end{Lemma}
\begin{Proof}
The proof is similar to that of Lemma~\ref{degwt}.
\end{Proof}

We are ready to state our main theorem.
\begin{Theorem} 
\label{mainthm}
\begin{enumerate}
\item The covariants form a free $K[c_4,c_6]$-module of rank~$2$ generated
by covariants $U$ and $H$ of weights $0$ and $2$. 
\item The contravariants form a free $K[c_4,c_6]$-module of rank~$2$ generated
by contravariants $P$ and $Q$ of weights $4$ and $6$. 
\end{enumerate}
\end{Theorem}

Our labelling of the covariants as $U$ and $H$, and
contravariants as $P$ and $Q$, follows the notation used
by Salmon \cite[Arts~217-221]{Salmon} in the case $n=3$.
The covariant $U$ is the identity map.
We call $H$ the Hessian since in degrees $2$ and $3$ it is
computed as the determinant of the matrix of second partial 
derivatives. We know of no such simple construction in degrees
$4$ and $5$. Since the Hessian of a genus one model 
is again a genus one model there is no natural generalisation of the 
statement (specific to the case $n=3$) that a plane cubic and its 
Hessian meet at the points of inflection of the cubic.

\section{The discrete invariants}
\label{sec:discinv}

We recall some classical theory 
from Klein's {\em Lectures on the
icosahedron} \cite{Kicosa}. 

\begin{Definition}
Let $\Delta_n$ be the subgroup of $\PGL_2$ generated by
$$ \begin{array}{llcl}  \smallskip
n=2 \qquad & 
\begin{pmatrix} 1 & 1/8 \\ 24 & -1 \end{pmatrix} & \text{and} &
\begin{pmatrix} 1 & 0 \\ 0 & -1 \end{pmatrix} \\ \smallskip
n=3 &
\begin{pmatrix} 1 & 1/3 \\ 6 & -1 \end{pmatrix} & \text{and} &
\begin{pmatrix} 1 & 0 \\ 0 & \zeta_3 \end{pmatrix} \\ \smallskip
n=4 &
\begin{pmatrix} 1 & 1/2 \\ 2 & -1 \end{pmatrix} & \text{and} &
\begin{pmatrix} 1 & 0 \\ 0 & \zeta_4 \end{pmatrix} \\
n=5 & 
\begin{pmatrix} \varphi & 1 \\ 1 & -\varphi \end{pmatrix} & \text{and} &  
\begin{pmatrix} 1 & 0 \\ 0 & \zeta_5 \end{pmatrix}
\end{array} $$
where $\zeta_n$ is a primitive $n$th root of unity 
and $\varphi = 1 + \zeta_5+\zeta_5^4$.
\end{Definition}

For $n=3,4,5$ the group $\Delta_n$ acts on $\PP^1$ as
the group of rotations of a tetrahedron, octahedron, icosahedron. 
Under stereographic projection the vertices of these Platonic solids 
are at the roots of
\begin{align*}
n=2 \qquad \quad D & = a(64 a^2-b^2) \\
n=3 \qquad \quad D & =  -a(27 a^3+b^3) \\
n=4 \qquad \quad D & =  ab(16a^4-b^4) \\
n=5 \qquad \quad D & =  ab(a^{10}-11 a^5 b^5-b^{10}). 
\end{align*}
The midpoints of the faces and edges are at the roots of
\begin{equation}
\label{c4def}
c_4  = 
\tfrac{-1}{((\deg D)-1)^2} \left| \begin{matrix} \smallskip
\frac{\partial^2 D}{\partial a^2}&\frac{\partial^2 D}{\partial a \partial b} \\
\frac{\partial^2 D}{\partial a \partial b}&\frac{\partial^2 D}{\partial b^2} 
\end{matrix} \right| 
\end{equation}
and
\begin{equation}
\label{c6def}
c_6  =  
\tfrac{1}{\deg c_4} \left| \begin{matrix} \smallskip
\frac{\partial D}{\partial a}&\frac{\partial D}{\partial b} \\
\frac{\partial c_4}{\partial a}&\frac{\partial c_4}{\partial b} 
\end{matrix} \right|. 
\end{equation}

\begin{Definition}
\label{defGamma}
Let $\GGamma_n$ be the inverse image of $\Delta_n$ in
$\SL_2$, and let $\Gamma_n$ be the commutator
subgroup of $\GGamma_n$.
The ring of {\em discrete invariants} is 
$$ K[a,b]^{\Gamma_n} = \{ f \in K[a,b] : f \circ \gamma = f 
\text{ for all } \gamma \in \Gamma_n \}. $$
\end{Definition}

\begin{Theorem}[Klein]
\label{Klein}
The ring of discrete invariants is generated by 
$c_4$, $c_6$ and $D$, subject only to
the relation $c_4^3-c_6^2 = 1728 D^n.$
\end{Theorem}
\begin{Proof} 
Since the characteristic of $K$ does not divide the order of $\Gamma_n$
this is a standard calculation.
We checked the answer using {\sf MAGMA} \cite{magma}.
\end{Proof}

We describe the action of $\GGamma_n$ on the discrete invariants.

\begin{Lemma}
\label{chilemma}
There is a unique character $\chi : \GGamma_n \to \Gm$
of order $6-n$ such that  
\begin{equation}
\label{gammac4c6} \begin{aligned}
D \circ \gamma & =  \chi(\gamma) D \\
c_4 \circ \gamma & =  \chi(\gamma)^2 c_4 \\
c_6 \circ \gamma & =  \chi(\gamma)^3 c_6 
\end{aligned} \end{equation}
for all $\gamma \in \GGamma_n$. Moreover $\ker(\chi) = \Gamma_n$. 
\end{Lemma}
\begin{Proof}
This follows by direct calculation.
\end{Proof}

\section{The Hesse family}
\label{sec:hesse}

In Section~\ref{sec:inv} we defined the ring of invariants $K[X_n]^{G_n}$ 
and in Section~\ref{sec:discinv} we defined the ring of 
discrete invariants $K[a,b]^{\Gamma_n}$.
We now identify $K[X_n]^{G_n}$ as a subring
of $K[a,b]^{\Gamma_n}$. To do this we 
first define a linear map $u_n : K^2 \to X_n$,
\begin{align*} \smallskip
u_2(a,b) & =  a(x^4+z^4) + b (\tfrac{1}{4} x^2 z^2) \\  \smallskip
u_3(a,b) & =  a(x^3+y^3+z^3) + b xyz \\ \smallskip
u_4(a,b) & =  \begin{pmatrix} a (x_1^2+x_3^2)- b x_2 x_4 \\
a (x_2^2+x_4^2)- b x_1 x_3 \end{pmatrix} \\ 
u_5(a,b) & =  \begin{pmatrix}
0 & a x_1 & b x_2 & -b x_3 & -a x_4 \\ 
 & 0 & a x_3 & b x_4 & -b x_5 \\ 
 & & 0 & a x_5 & b x_1 \\ 
 & - & & 0 & a x_2 \\ 
 &   & & & 0 
\end{pmatrix}.  
\end{align*}
The models $u_n(a,b)$ are called
{\em Hesse models}. Collectively they form the {\em Hesse family}.
The geometry of the Hesse family is discussed, for example,
in \cite{ArtDol}, \cite{BHM} in the cases $n=3,5$.
The following two propositions will be proved in Section~\ref{sec:heis}.

\begin{Proposition}
\label{hcovj}
Every non-singular model $\phi \in X_n$ is equivalent to
a Hesse model.
\end{Proposition}

\begin{Proposition}
\label{hlem}
Let $\G_n$ and $\GGamma_n$ be the groups defined in Sections~\ref{sec:inv}
and~\ref{sec:discinv}. 
\begin{enumerate}
\item There exists $g \in \G_n$ with $g \circ u_n = u_n$ 
and $\det(g)=-1$. 
\item For each $\gamma \in \GGamma_n$ there exists $g \in \G_n$
with $g \circ u_n = u_n \circ \gamma$. 
\end{enumerate}
\end{Proposition}

The map $u_n : K^2 \to X_n$ induces a homomorphism
of polynomial rings $u_n^*: K[X_n] \to K[a,b] \, ; \,
 F \mapsto F \circ u_n$. Proposition~\ref{hlem}(ii)
implies an analogous result where $\G_n$ and $\GGamma_n$ are replaced
by their commutator subgroups $G_n$ and $\Gamma_n$.
It follows that $u_n^*$  restricts to a map
\[ u_n^* : K[X_n]^{G_n} \to K[a,b]^{\Gamma_n}. \]

\begin{Lemma} 
\label{INJ}
The map $u_n^* : K[X_n]^{G_n} \to K[a,b]^{\Gamma_n}$
is injective.
\end{Lemma}
\begin{Proof} 
Let $F \in K[X_n]^{G_n}$ be a homogeneous invariant vanishing on the 
Hesse family. By Proposition~\ref{hcovj} 
it also vanishes at every non-singular $\phi \in X_n$.
By Theorem~\ref{invthm}(ii) the latter are Zariski dense
in $X_n$. It follows that $F$ is identically zero.
\end{Proof}

\begin{Lemma}
\label{theymatch}
The map $u_n^*$ takes the invariants $c_4$, $c_6$ and $\Delta$ of
Theorem~\ref{invthm} to the discrete invariants $c_4$, $c_6$ and
$D^n$ of Theorem~\ref{Klein}.
\end{Lemma}
\begin{Proof}
We compute the invariants of the generic Hesse model using the formulae
and algorithms in \cite[Sections~7 and~8]{g1inv}. This gives an alternative
computational proof of Lemma~\ref{INJ}.
\end{Proof}

\begin{Lemma}
\label{chidet}
If $g \in \G_n$ and $\gamma \in \GGamma_n$ satisfy  
$g \circ u_n = u_n \circ \gamma$ then $\chi(\gamma) = (\det g)^2$.
\end{Lemma}
\begin{Proof} The map $u_n^*$ identifies the invariants 
$c_4$ and $c_6$ with the corresponding discrete invariants.
Since the former have weights $4$ and $6$, and the latter
satisfy~(\ref{gammac4c6}), we deduce
\begin{align*}
\chi(\gamma)^2 c_4(a,b) & =  (\det g)^4 c_4(a,b) \\
\chi(\gamma)^3 c_6(a,b) & =  (\det g)^6 c_6(a,b). 
\end{align*}
It follows that $\chi(\gamma) = (\det g)^2$.
\end{Proof}
 
We say that a discrete invariant is an invariant if it belongs
to the image of $u_n^*$. The following theorem characterises the
invariants among the discrete invariants, and thus serves as a 
prototype for our treatment of the covariants and contravariants
in Section~\ref{sec:covcontra}. 

\begin{Theorem}
\label{proto}
Let $f$ be a homogeneous discrete invariant of degree $d$. 
Then $f$ is an invariant if and only if $d=kn/(6-n)$ for
some even integer $k$ and 
\begin{equation}
\label{fpropinv}
f \circ \gamma = \chi (\gamma)^{k/2} f 
\end{equation}
for all $\gamma \in \GGamma_n$.
\end{Theorem}

\begin{Proof}
Suppose that $f = F \circ u_n$ for some invariant $F$.
Lemma~\ref{degwt} shows that since $F$ is homogeneous of 
degree $d$ it has weight $k$ where $d=kn/(6-n)$. 
We use Proposition~\ref{hlem}(i) to show that 
$k$ is even, and then Proposition~\ref{hlem}(ii) combined with
Lemma~\ref{chidet} to establish~(\ref{fpropinv}).

For the converse we use the description of 
the discrete invariants given in Theorem~\ref{Klein},
namely that $K[a,b]^{\Gamma_n}$ is a 
free $K[c_4,c_6]$-module of rank $n$ with basis $1,D, \ldots, D^{n-1}$.
Using Lemma~\ref{chilemma} 
we find that $c_4$ and $c_6$ satisfy the conditions 
required of an invariant, but $D$, $D^2$, $\ldots$, $D^{n-1}$ do not.
It only remains to show that there are invariants of
weights $4$ and $6$. This was established in Theorem~\ref{invthm}.
\end{Proof}

\begin{Remark}
Our use of the Hesse family in the above proof is analogous 
to our use of the Weierstrass family in the proof of 
Theorem~\ref{invthm}. The advantage of the Weierstrass family
is that it allows us to work without restriction on the characteristic of $K$.
The advantage of the Hesse family is that it allows us 
to study the covariants and contravariants.
\end{Remark}

\section{The discrete covariants}
\label{sec:disccov}

In Section~\ref{sec:discinv} we defined subgroups 
$\Gamma_n \subset \SL_2$ for $n=2,3,4,5$.

\begin{Definition} A {\em discrete covariant} is a $\Gamma_n$-equivariant
polynomial map $p:K^2 \to K^2$. It is represented by a pair of 
polynomials $(p_1,p_2)$ with  $p_1,p_2 \in K[a,b]$.
\end{Definition}

The discrete covariants form a module $M$ over the ring of 
discrete invariants $R=K[a,b]^{\Gamma_n} = K[D,c_4,c_6]$.
There is a derivation 
$$ \begin{array}{ll} \partial : R \to M \,; & f
 \mapsto (-\frac{\partial f}{\partial b}, \frac{\partial f}{\partial a})
\end{array} $$
and an $R$-bilinear alternating form 
$$ \begin{array}{ll} [~,~] : M \times M \to R \,; & 
(p,q) \mapsto p_1 q_2 -p_2 q_1. \end{array}$$
We write $U=(a,b)$ for the identity map.

\begin{Theorem}
\label{thm:disccov}
The discrete covariants form a free $K[c_4,c_6]$-module of rank $2n$ with 
basis $D^i U$, $D^i \partial D$, 
$\partial c_4$, $\partial c_6$ for $i=0,1, \ldots ,n-2$.
\end{Theorem}

For the proof we first show that $M$ is a free $K[c_4,c_6]$-module.
Then we show, by computing the Hilbert series, that the elements listed
have the right degrees to be generators. Finally we check that
our putative basis is independent.

\begin{Lemma} $M$ is a free $K[c_4,c_6]$-module.
\end{Lemma}
\begin{Proof}
The proof follows the method described in \cite[Section~4.3]{Benson}.

Since $c_4$ and $c_6$ are coprime they form a regular sequence
in $K[a,b]$, and so $K[a,b]$ is a free $K[c_4,c_6]$-module.
The projection map
$$ K[a,b]^2 \to M \, ; \quad  p \mapsto 
\frac{1}{|\Gamma_n|} \sum_{\gamma \in \Gamma_n} 
\gamma \circ p \circ \gamma^{-1} $$
shows that $M$ is a projective $K[c_4,c_6]$-module. 
By \cite[Theorem 4.1.1]{Benson} it is therefore a free
$K[c_4,c_6]$-module. 
\end{Proof}

\begin{Lemma}
\label{poincare}
The Hilbert series of $M$ is
$$ h_M(z) = \frac{(z^{r-1}+z^{s-1}) + \sum_{i=0}^{n-2} (z^{it+1}+z^{it+t-1}) 
}{(1-z^r)(1-z^s)} $$
where 
$r = \deg c_4$, 
$s = \deg c_6$  
and $t=\deg D$. 
\end{Lemma}
\begin{Proof}
By Molien's theorem \cite[Theorem 2.5.3]{Benson} 
the Hilbert series of $R$ and $M$ are 
$$ h_R (z) 
= \frac{1}{|\Gamma_n|}
\sum_{\gamma \in \Gamma_n} \frac{1}{1- \Tr (\gamma) z + z^2} $$
and 
$$ h_M (z) 
=  \frac{1}{|\Gamma_n|}
\sum_{\gamma \in \Gamma_n} \frac{\Tr(\gamma)}{1- \Tr (\gamma) z + z^2}. $$
Thus
\begin{equation*}
\label{poincare1}
(1+z^2)h_R(z) = 1+ z h_M(z).
\end{equation*}
But by Theorem~\ref{Klein} we already have
\begin{equation*}
\label{poincare2}
h_R(z) = \frac{\sum_{i=0}^{n-1} z^{it}}{(1-z^r)(1-z^s)}.
\end{equation*}
The lemma follows on noting that $(n-1)t+1= r+s-1$. (In fact 
$r=4n/(6-n)$, $s= 6n/(6-n)$ and $t=12/(6-n)$.)
\end{Proof}

By Lemma~\ref{poincare} the discrete covariants
listed in the statement of Theorem~\ref{thm:disccov} have the right 
degrees to generate $M$ as a $K[c_4,c_6]$-module. 
It remains to show that they are independent.

\begin{Lemma}
\label{M'}
Let $f$ be a non-zero homogeneous discrete invariant. 
If $\charic(K) \nmid \deg (f)$ then the discrete covariants 
$U$ and $\partial f$ generate a free $R$-module $M_f$ 
with $fM \subset M_f \subset M$.
\end{Lemma}
\begin{Proof}
The module $M_f$ is free since by Euler's identity we have
$$[U, \partial f] = (\deg f) f \not= 0.$$
It contains $fM$ since for $p \in M$ we have
\begin{equation}
\label{eq0}
 (\deg f) f p = -[\partial f,p] U + [ U,p ] \partial f .
\end{equation}
\end{Proof}
In the notation of this section we may 
re-write~(\ref{c6def})  
as
\begin{equation}
\label{eq1}
 [\partial D ,\partial c_4] = (\deg c_4) c_6. 
\end{equation}
Applying $\partial$ and then $[\partial D,-]$ to the relation of
Theorem~\ref{Klein} we obtain
\begin{equation}
\label{eq2}
 [\partial D ,\partial c_6] = (\deg c_6) c_4^2. 
\end{equation}
Then we take $f=D$ in~(\ref{eq0}) to get
\begin{equation}
\label{eq3}
\begin{aligned}
3 D \partial c_4 & =  n(-c_6 U + c_4 \partial D) \\
2 D \partial c_6 & =  n(-c_4^2 U + c_6 \partial D). 
\end{aligned} 
\end{equation}
Taking linear combinations and 
using the identity $c_4^3-c_6^2 = 1728 D^n$ we find
\begin{equation}
\label{eq4}
\begin{aligned}
1728 n D^{n-1} U & =  3 c_6 \partial c_4 - 2 c_4 \partial c_6 \\
1728 n D^{n-1} \partial D & =  3 c_4^2 \partial c_4 - 2 c_6 \partial c_6. 
\end{aligned}
\end{equation}

We already know by Theorem~\ref{Klein} that $R$ is a free 
$K[c_4,c_6]$-module of rank $n$, with basis $1, D, \ldots ,D^{n-1}$. 
Taking $f=D$ in 
Lemma~\ref{M'} it follows that $M$ is a $K[c_4,c_6]$-module 
of rank $2n$. Moreover the discrete covariants
$D^i U, D^i \partial D$ for $i=0,1, \ldots, n-1$ generate 
a free submodule of maximal rank. The relations
(\ref{eq3}) and~(\ref{eq4}) show that we can replace $D^{n-1} U$ and 
$D^{n-1} \partial D$ by $\partial c_4$ and $\partial c_6$ without 
destroying this property. In other words the discrete 
covariants listed in the statement of Theorem~\ref{thm:disccov}
are independent. 
This completes the proof of Theorem~\ref{thm:disccov}.

\medskip

We describe the action of $\GGamma_n$ on the 
discrete covariants. The character $\chi:\GGamma_n \to \Gm$ 
was defined in Lemma~\ref{chilemma}.
\begin{Lemma}
\label{chicov}
If a discrete invariant $f$ satisfies
$$ f \circ \gamma = \chi(\gamma)^r f $$
for all $\gamma \in \Gamma_n$,
then the discrete covariant $p=\partial f$ satisfies
$$ p \circ \gamma = \chi(\gamma)^r \gamma  \circ p $$
for all $\gamma \in \Gamma_n.$
\end{Lemma}
\begin{Proof}
This is proved by a short calculation using the chain rule.
\end{Proof}

\section{The covariants and contravariants}
\label{sec:covcontra}


We use the discrete covariants to study
the covariants and contravariants (see Definition~\ref{def:cov}), 
just as in Section~\ref{sec:hesse}
we used the discrete invariants to study the invariants.
Recall that for $g \in \G_n$ we write $g^T$ for the element obtained by
transposing all constituent matrices. 
The following proposition will be proved in Section~\ref{sec:heis}.

\begin{Proposition}
\label{hn1}
There is a finite subgroup $H_n \subset G_n$ with the following properties
\begin{enumerate}
\item $H_n = \{ g \in G_n : g \circ u_n = u_n \}$ 
\item The image of $u_n$ is  $X_n^{H_n} = \{ \phi \in X_n : 
g \phi = \phi \text{ for all } g \in H_n \}.$ 
\item If $g \in H_n$ then $g^T \in H_n$.
\end{enumerate}
\end{Proposition}

It is possible to view the contravariants as $G_n$-equivariant
polynomial maps from $X_n$ to its dual $X_n^*$. The connection is 
afforded by the following pairing on $X_n$.

\begin{Lemma}
\label{pair}
There is a symmetric bilinear form $\langle~,~\rangle$ on $X_n$ such that \\
(i) $ \langle g \phi_1, \phi_2 \rangle = \langle \phi_1, g^T \phi_2 \rangle $
for all $g \in \G_n$ and $\phi_1,\phi_2 \in X_n$, \\
(ii) $\langle~,~\rangle$ is non-degenerate on the image of $u_n$.
\end{Lemma}
\begin{Proof}
A suitable pairing is
\begin{align*} \smallskip
n=2 \qquad  \langle f, g \rangle & =  
f(\tfrac{\partial}{\partial x},\tfrac{\partial}{\partial z}) 
g(x,z) \\ \smallskip
n=3 \qquad \langle f, g \rangle & =  
f(\tfrac{\partial}{\partial x},\tfrac{\partial}{\partial y},
  \tfrac{\partial}{\partial z}) 
g(x,y,z) \\  \smallskip
n=4 \qquad \langle f, g \rangle 
& =   
\textstyle\sum_{i=1}^2 
f_i(\tfrac{\partial}{\partial x_1}, \ldots ,
  \tfrac{\partial}{\partial x_4})  g_i(x_1, \ldots ,x_4) \\
n=5 \qquad  \langle f, g \rangle 
& =   \textstyle\sum_{i<j} 
f_{ij}(\tfrac{\partial}{\partial x_1}, \ldots,
\tfrac{\partial}{\partial x_5}) 
g_{ij}(x_1, \ldots ,x_5). 
\end{align*}
Properties (i) and (ii) are checked by routine calculation.
\end{Proof}

The Hesse family was defined in Section~\ref{sec:hesse} by specifying
a map $u_n : K^2 \to X_n$. We now view $K^2$ as a space of column vectors.
By Lemma~\ref{pair}(ii) there is a matrix $\eps_n \in \GL_2(K)$
such that
\begin{equation}
\label{epsdef}
 \langle u_n(x), u_n(\eps_n y) \rangle = x_1 y_2 - x_2 y_1 
\end{equation}
for all $x,y \in K^2$. The following lemma is 
required for our treatment of the contravariants. 
\begin{Lemma}
\label{epslem}
If $g \in \G_n$ and $\gamma \in \GGamma_n$ satisfy
$ g \circ u_n = u_n \circ \gamma $
then
$$ g^{-T} \circ u_n \circ \eps_n = u_n \circ \eps_n \circ \gamma. $$
\end{Lemma}
\begin{Proof} 
For $x,y \in K^2$ we have 
\begin{align*}
\langle u_n(x),u_n(\eps_n y) \rangle 
& =  \langle u_n(\gamma x ),u_n( \eps_n \gamma y)  \rangle
 && \text{ by~(\ref{epsdef}) and $\GGamma_n \subset \SL_2$}  \\
& =  \langle g (u_n (x)),u_n (\eps_n \gamma y) \rangle 
&& \text{ since $g \circ u_n = u_n \circ \gamma$ } \\
& =  \langle u_n(x),g^T (u_n (\eps_n \gamma y)) \rangle
&& \text{ by Lemma~\ref{pair}(i). }
\end{align*}
Since $g$ acts on the image of $u_n$, 
Proposition~\ref{hn1}(i) gives $g H_n g^{-1} = H_n$.
We deduce by Proposition~\ref{hn1}(iii) that $g^T H_n g^{-T} = H_n$
and hence by Proposition~\ref{hn1}(ii) that $g^T$ acts on the 
image of $u_n$. The lemma now follows by Lemma~\ref{pair}(ii)
and the above calculation.
\end{Proof}

\begin{Proposition}
\label{covINJdisccov}
Let $F: X_n \to X_n$ be a covariant, respectively contravariant. 
Then there is a discrete covariant $f$ such that 
$F \circ u_n  = u_n  \circ f$, 
respectively $F \circ u_n  = u_n  \circ \varepsilon_n \circ f$.
Moreover $F$ is uniquely determined by $f$. 
\end{Proposition}
\begin{Proof}
Proposition~\ref{hn1} shows that $F$ acts on the image of $u_n$.
So there is a polynomial map $f: K^2 \to K^2$ satisfying
$F \circ u_n  = u_n  \circ f$, respectively
$F \circ u_n  = u_n  \circ \varepsilon_n \circ f$.
It follows by Proposition~\ref{hlem}(ii), combined with Lemma~\ref{epslem}
in the case $F$ is a contravariant, that $f$ is a discrete covariant.

If $F_1$ and $F_2$ determine the same discrete covariant $f$ 
then by Proposition~\ref{hcovj} they agree on all non-singular models.
By Theorem~\ref{invthm}(ii) the non-singular models are Zariski dense
in $X_n$, and from this we deduce that $F_1= F_2$.
\end{Proof}
We say that a discrete covariant $f$ is a covariant, respectively contravariant, 
if it arises as described in Proposition~\ref{covINJdisccov}.
We obtain the following analogue of Theorem~\ref{proto}.
\begin{Theorem} 
\label{describecov}
Let $f$ be a homogeneous 
discrete covariant of degree $d$. Then $f$ is a covariant,
respectively contravariant, 
if and only if $d=1+kn/(6-n)$, respectively $d=-1+kn/(6-n)$, 
for some even integer $k$ and 
\begin{equation}
\label{fprop1}
f \circ \gamma = \chi (\gamma)^{k/2} \gamma \circ f 
\end{equation}
for all $\gamma \in \GGamma_n$. 
\end{Theorem}
\begin{Proof}
Suppose that $F \circ u_n = u_n \circ f$, respectively 
$F \circ u_n = u_n \circ \eps_n \circ f$, for some 
covariant, respectively contravariant, $F$.  
Lemma~\ref{prop:covwts} shows that, since $F$ is homogeneous of 
degree $d$, it has weight $k$ 
where $d=1+kn/(6-n)$, respectively $d=-1+kn/(6-n)$. 
In the case $F$ is a covariant
we use Proposition~\ref{hlem}(i) to show that $k$ is 
even, and Proposition~\ref{hlem}(ii) combined with Lemma~\ref{chidet}
to establish~(\ref{fprop1}). In the case $F$ is a contravariant
we use Lemma~\ref{epslem} to make the necessary modifications.

For the converse we use the description of the discrete covariants 
given in Theorem~\ref{thm:disccov}, namely that $M$ is 
a free $K[c_4,c_6]$-module generated by $D^i U$, $D^i \partial D$,
$\partial c_4$, $\partial c_6$ 
for $i=0, \ldots, n-2$. 
Using Lemmas~\ref{chilemma} and~\ref{chicov} we find that 
only $U$ and $\partial D$ satisfy the conditions
required of a covariant (with $k=0,2$) and only 
$\partial c_4$ and $\partial c_6$ satisfy the conditions
required of a contravariant (with $k=4,6$).
To complete the proofs of Theorems~\ref{mainthm} and~\ref{describecov},  
it remains to show that there are covariants of weights $0$
and $2$ and contravariants of weights $4$ and $6$.

We construct the contravariants using a method described
by Salmon in the case $n=3$; see \cite[Arts.~220, 221]{Salmon}. 
By Lemma~\ref{pair} we may identify the
contravariants with the space $\Pol_{G_n}(X_n,X_n^*)$
of $G_n$-equivariant polynomial maps from $X_n$ to its dual $X_n^*$.
If we pick a basis $x_1, \ldots, x_N$ for $X_n^*$ then 
$K[X_n] = K[x_1, \ldots, x_N]$ and there is a derivation
$$ \delta : K[X_n]^{G_n} \to \Pol_{G_n}(X_n,X_n^*) $$
given by $(\delta F)(\phi) = \sum_{i=1}^N \frac{\partial F}{\partial x_i}
(\phi) x_i$. It may be checked that $\delta$ is independent of
the choice of basis $x_1, \ldots, x_N$.
The contravariants of weights $4$ and $6$ are the so-called evectants
$\delta c_4$ and $\delta c_6$ of the invariants $c_4$ and
$c_6$ in Theorem~\ref{invthm}.

The covariant of weight $0$ is of course the identity map. 
So it only remains to show that there is a covariant of weight $2$.
It is clear from Definition~\ref{def:cov} that the composition
of two contravariants is a covariant. If $n=2$ then 
the contravariants have degrees $1$ and $2$, and their composition
is the required covariant of weight $2$. Otherwise, composing the
contravariant of weight $4$ with itself we find that
$$ f_n(c_4,c_6) U + g_n(c_4,c_6) \partial D$$
is a covariant where
\begin{align*}
f_3(c_4,c_6) & =  3 c_4^2 
&  f_4(c_4,c_6) & =  81 c_4^6 + 40 c_4^3 c_6^2 - c_6^4 \\
g_3(c_4,c_6) & =  c_6  &  
g_4(c_4,c_6) & =  6 c_4 c_6 (5 c_4^3-c_6^2) 
\end{align*} 
and 
\begin{align*}
f_5(c_4,c_6)  = & 184528125 c_4^{18} + 230364000 c_4^{15} c_6^{2} - 25697763 c_4^{12} c_6^{4} \\ & ~\qquad\qquad\medskip~+4909960 c_4^{9} c_6^{6} 
+ 44583 c_4^{6} c_6^{8} + 984 c_4^{3} c_6^{10} - c_6^{12} \\
g_5(c_4,c_6)  = & 18 c_4 c_6 ( 2399625 c_4^{15} - 658917 c_4^{12} c_6^2  \\
& ~\qquad\qquad\medskip~+245498 c_4^9 c_6^4 + 4246 c_4^6 c_6^6 +
    205 c_4^3 c_6^8 - c_6^{10}). 
\end{align*}
Since $U$ is a covariant it follows that $g_n(c_4,c_6) \partial D$
is a covariant.

Let $f$ be a homogeneous discrete invariant of positive degree with
$D \nmid f$. We claim that if $f \partial D$ is a covariant
then $f_1 \partial D$ is also a covariant for some proper factor $f_1$
of $f$. To see this let
$(a:b)$ be a root of $f$. 
Then $\phi = u_n(a,b)$ is non-singular since $D(a,b) \not= 0$.
By \cite[Lemma 4.10]{g1inv} the Zariski closure of the orbit of $\phi$ 
is the zero locus of an irreducible invariant $F$. 
The covariant corresponding to $f \partial D$ vanishes on
the orbit of $\phi$ and is therefore divisible by $F$.
This proves the claim.

Finally we check for $n=3,4,5$ that $g_n(c_4,c_6)$ is not 
divisible by $\Delta = (c_4^3-c_6^2)/1728$, 
equivalently $g_n(1,1) \not=0$. 
In fact $g_3(1,1)=1$, $g_4(1,1) = 2^3 \cdot 3$ and $g_5(1,1) = 2^{14} \cdot
3^7$. It follows by the claim in the last paragraph 
that $\partial D$ is a covariant.
\end{Proof}

We write $U: X_n \to X_n$ for the identity map. 

\begin{Definition}
(i) The Hessian $H : X_n \to X_n$ is the unique covariant
(of weight 2) satisfying
$$ H \circ u_n = u_n \circ \partial D. $$ 
(ii) The contravariants $P,Q : X_n \to X_n$ are the unique
contravariants (of weights 4 and 6) satisfying
\begin{align*}
\kappa^{-1} (\deg c_4) P \circ u_n & =  u_n \circ \eps_n \circ \partial c_4 \\
\kappa^{-1} (\deg c_6) Q \circ u_n & =  u_n \circ \eps_n \circ \partial c_6. 
\end{align*}
where $\kappa = 1/4,1,2,5$ for $n=2,3,4,5$. (The scaling factor $\kappa$
has been chosen to simplify the formulae in Section~\ref{sec:formulae}.) 
\end{Definition}

\begin{Theorem}
\label{kappa}
Let $\langle~,~\rangle$ be the pairing defined in the proof of 
Lemma~\ref{pair}. Then
\begin{align*}
\langle U , P \rangle & =  \kappa c_4 &
\langle H , P \rangle & =  \kappa c_6 \\
\langle U , Q \rangle & =  \kappa c_6 & 
\langle H , Q \rangle & =  \kappa c_4^2. 
\end{align*}
\end{Theorem}
\begin{Proof}
It suffices to prove this for $\phi \in X_n$ a Hesse model. 
By~(\ref{epsdef}) the required identities are
\begin{align*}
[ U , \partial c_4 ] & =  (\deg c_4) c_4 & 
[ \partial D , \partial c_4 ] & =  (\deg c_4) c_6 \\
[ U , \partial c_6 ] & =  (\deg c_6) c_6 & 
[ \partial D , \partial c_6 ] & =  (\deg c_6) c_4^2. 
\end{align*}
These were proved in Section~\ref{sec:disccov}. 
\end{Proof}

\section{The Heisenberg group}
\label{sec:heis}

In this section we prove some results postponed
from Sections~\ref{sec:hesse} and~\ref{sec:covcontra}.
\begin{Definition}
A genus one normal curve $C \to \PP^{n-1}$ is \\
(i) if $n=2$ a double cover of $\PP^1$ ramified at 4 points, \\
(ii) if $n \ge 3$ a genus one curve embedded in $\PP^{n-1}$ 
by a complete linear system of degree $n$.
\end{Definition}

It is well known that $E =\Jac(C)$ acts on $C$ by translation, and 
translation by $P \in E$ extends to an automorphism of $\PP^{n-1}$
if and only if $P \in E[n]$. 

\begin{Definition} 
\label{heisdef}
The Heisenberg group of $C \to \PP^{n-1}$ is the group of
all matrices in $\SL_n$ that act on $C$ as translation
by an $n$-torsion point of its Jacobian.
As a group it is a central extension of $E[n]$ by $\mu_n$ with
commutator given by the Weil pairing
$e_n : E[n] \times E[n] \to \mu_n$.
\end{Definition}

\begin{Definition} The standard Heisenberg group of degree $n$
is the subgroup $H_n \subset \SL_n$ generated by
$$ \sigma_n = \xi_n \begin{pmatrix}
    1   &    0   &    0    & \cdots &      0     \\
    0   & \zeta_n  &    0    & \cdots &      0     \\
    0   &    0   & \zeta_n^2 & \cdots &      0     \\ 
 \vdots & \vdots & \vdots  &        &   \vdots   \\
    0   &   0    &    0    & \cdots & \zeta_n^{n-1}  
\end{pmatrix}
, \text{ and  }
\tau_n = \xi_n \begin{pmatrix}
    0   &    0   & \cdots &   0    &      1     \\ 
    1   &    0   & \cdots &   0    &      0     \\
    0   &    1   & \cdots &   0    &      0     \\ 
 \vdots & \vdots &        & \vdots &   \vdots   \\
    0   &   0    & \cdots &   1    &      0       
\end{pmatrix}, $$
where $\xi_n =1$ for $n$ odd, and $\xi_n = \zeta_{2n}$ for $n$ even.
\end{Definition}

\begin{Lemma} 
\label{choosecoords}
Let $C \to \PP^{n-1}$ be a genus one normal curve
with Jacobian $E$. Let $S$, $T$ be a basis for $E[n]$ with 
$e_n(S,T)=\zeta_n$. Then we can change co-ordinates on $\PP^{n-1}$ 
so that translation by $S$ and $T$ is given by the images
of $\sigma_n$ and $\tau_n$ in $\PGL_n$. In particular
$C \to \PP^{n-1}$ has Heisenberg group $H_n$.
\end{Lemma}
\begin{Proof}
This is standard. See for example \cite[Proposition 2.3]{jems}.
\end{Proof}


The following lemma is based on results in \cite[Chapter III]{Hu}.

\begin{Lemma}
\label{heis2}
Let $C \to \PP^{n-1}$ be a genus one normal curve
with Heisenberg group $H_n$. If $n=2,3,4,5$ then $C = C_\phi$ 
for some Hesse model~$\phi$. 
\end{Lemma}

\begin{Proof} Case $n=2$. We decompose $X_2$ as an $H_2$-module
and find that (up to scalars) there are exactly three binary quartics 
whose roots are permuted by $H_2$, but do not belong to the Hesse family.
These are the pairwise products of $x z$, $x^2- z^2$ 
and $x^2+z^2$. These quartics do not have Heisenberg group $H_2$,
since $H_2$ fails to act transitively on their roots.

Case $n=3$. We decompose $X_3$ as a $H_3$-module and find
that (up to scalars) there are exactly eight ternary cubics
that define curves fixed by $H_3$, but do not belong to the 
Hesse family. These are 
\[ \begin{array}{ll}
x^3 + \zeta_3^i y^3 + \zeta_3^{2 i} z^3 & \text{ for } i =1,2 \\
x^2 y + \zeta_3^i y^2 z + \zeta_3^{2 i} x z^2 & 
 \text{ for } i =0,1,2 \\
x y^2 + \zeta_3^i y z^2 + \zeta_3^{2 i} x^2 z & 
 \text{ for } i =0,1,2.
\end{array} \]
These curves do not have Heisenberg group $H_3$ 
since the action of $H_3$ modulo its centre is not fixed point free.

Case $n=4$. We decompose the space of quadrics in $4$ variables as
an $H_4$-module. We find that there are exactly three $2$-dimensional
subspaces that define a curve fixed by $H_4$, but do not
belong to the Hesse family. These are spanned by 
\[ \begin{array}{lcl}
x_1 x_2 + x_3 x_4 & \text{and} &  x_2 x_3 + x_1 x_4,  \\
x_1 x_2 - x_3 x_4  & \text{and} & x_2 x_3 - x_1 x_4,  \\
x_1^2 - x_3^2  & \text{and} & x_2^2 - x_4^2 .
\end{array} \]
Each of these pairs of quadrics defines a singular curve.

Case $n=5$. We take $C \subset \PP^4$ with equations
$$ a x_i^2 + b x_{i+1} x_{i+4} + c x_{i+2} x_{i+3} =0 
\quad \text{ for } i=1,2,3,4,5 $$
where all subscripts are read mod 5. Since $C$ has Heisenberg
group $H_5$ it meets the hyperplane $\{x_1=0\}$ in 5 distinct points.
By Riemann-Roch these points span the hyperplane. So  
$C$ contains a point of the form $(0:z_2:z_3:z_4:z_5)$
with each $z_i$ non-zero. A short calculation then shows that
$a^2 + b c = 0$ and so $C=C_\phi$ where $\phi = u_5(a,b)$.
\end{Proof}

\begin{Lemma}
\label{existsA}
Let $\phi,\phi' \in X_n$ be non-singular models. If
$C_{\phi} = C_{\phi'}$ then $\phi$ and $\phi'$ are equivalent. 
Moreover if $n=4,5$ and $\phi' = [A,I_n] \phi$ 
then $A \in \GL_2$ is uniquely determined if $n=4$, and $A \in \GL_5$
is uniquely determined up to sign if $n=5$.
\end{Lemma}
\begin{Proof}
This is clear for $n=2,3,4$. The case $n=5$ follows from the 
Buchsbaum-Eisenbud acyclicity criterion and the properties of 
minimal free resolutions. See for example \cite[Section 5.2]{g1inv}.
\end{Proof}

Combining the last three lemmas shows that 
every non-singular model is equivalent to a Hesse model.

\medskip

\begin{ProofOf}{Proposition~\ref{hcovj}} 
Let $\phi \in X_n$ be a non-singular model. By definition this means that
$C_\phi$ is a smooth curve of genus one. In the cases $n=2,3$ it is clear
that $C_\phi \to \PP^{n-1}$ is a genus one normal curve. The cases
$n=4,5$ are treated in \cite[Proposition 5.10(i)]{g1inv}. 
By Lemma~\ref{choosecoords} we may assume that $C_\phi \to \PP^{n-1}$ has 
Heisenberg group $H_n$. Then Lemma~\ref{heis2} shows that 
$C_{\phi} = C_{\phi'}$ for some Hesse model $\phi'$ and finally 
Lemma~\ref{existsA} shows that $\phi$ and $\phi'$ are equivalent.
\end{ProofOf}

If $n=2,3$ then $H_n$ is already a subgroup of $G_n= \SL_n$.
If $n=4,5$ we identify $H_n$ as a subgroup of $G_n$ via
$\sigma_4 \mapsto [\sigma_2,\sigma_4]$, $\tau_4 \mapsto [\tau_2,\tau_4]$
and $\sigma_5 \mapsto [\zeta_5^3 \sigma_5^4 ,\sigma_5]$, 
$\tau_5 \mapsto [\tau_5^3,\tau_5]$.
\begin{Lemma}
\label{hcalc} 
Let $u_n : K^2 \to X_n$ be the linear map defining the
Hesse family. Then $g \circ u_n = u_n$ for all $g \in H_n$.
\end{Lemma}
\begin{Proof}
This is checked by direct calculation.
\end{Proof}

\begin{Lemma}
\label{heiscond}
Let $\phi \in X_n$ be a non-singular Hesse model. Then 
$C_\phi$ has Heisenberg group $H_n$ and 
$H_n = \{ g \in G_n : g \phi = \phi \}$.
\end{Lemma}
\begin{Proof}
In view of the last two lemmas, it suffices to show that
if $g \in G_n$ with $g \phi  =\phi$ then the automorphism
$\gamma$ of $C_{\phi}$ induced by $g$ is a translation map.
By \cite[Proposition 5.19]{g1inv} we have 
$\gamma^* \omega_{\phi}= \omega_\phi$. So this follows from 
\cite[Lemma 2.4]{g1inv}. 
\end{Proof}

Next we show that $H_n \subset G_n$ has the properties
stated in Section~\ref{sec:covcontra}.

\medskip

\begin{ProofOf}{Proposition~\ref{hn1}} 
(i) We must show that $H_n = \{g \in G_n : g \circ u_n = u_n \}$.
This follows from Lemmas~\ref{hcalc} and~\ref{heiscond}. \\
(ii) By Lemma~\ref{hcalc} we have $\im (u_n) \subset X_n^{H_n}$. 
We prove equality by showing that $\dim(X_n^{H_n})=2$.
The character of $X_n$ as a representation of $H_n$
is constant on the centre of $H_n$ and elsewhere
takes value $\xi_n=1,1,0$ for $n=2,3,5$. Thus 
$$\dim ( X_n^{H_n}) = \frac{1}{n^3} ( n \dim X_n + (n^3-n) \xi_n) = 2.$$
The case $n=4$ is similar. \\
(iii) It is clear from the definition of $H_n$ 
that if $g \in H_n$ then $g^T \in H_n$.
\end{ProofOf}

We prepare for the proof of Proposition~\ref{hlem} by describing
the normaliser of $H_n$, where $H_n$ is viewed first as a subgroup of 
$\GL_n$ and then as a subgroup of $\G_n$. 

\begin{Lemma}
\label{normaliser1}
There is an exact sequence 
$$ \xymatrix{ 0 \ar[r] & \Theta_n \ar[r] & N_{\GL_n}(H_n) \ar[r]^-{\pi} &
\SL_2(\Z/n\Z) \ar[r] & 0. } $$
where $\Theta_n \subset \GL_n$ is 
generated by $H_n$ and the scalar matrices.
\end{Lemma}
\begin{Proof}
Let $g \in \GL_n$ with $g H_n g^{-1} = H_n$. 
Writing $\propto$ for equality in $\PGL_n$ we have
$g \, \sigma_n \, g^{-1}  \propto  \sigma_n^a \, \tau_n^c$ and
$ g \, \tau_n \, g^{-1}  \propto  \sigma_n^b \, \tau_n^d$
for some $a,b,c,d \in \Z/n\Z$.
We define $\pi(g) = \left(
\begin{smallmatrix} a & b \\ c & d \end{smallmatrix}\right)$.
It is easy to check that $\pi$ is a group homomorphism with
kernel $\Theta_n$. Then Lemma~\ref{choosecoords} shows that 
$\im(\pi) = \SL_2(\Z/n\Z)$.
\end{Proof}

\begin{Lemma}
\label{pm1}
Let $C \subset \PP^{n-1}$ be a genus one normal curve
with Heisenberg group $H_n$ and $j(C) \not= 0,1728$.
If $g \in N_{\GL_n}(H_n)$ acts on $C$ then $\pi(g) = \pm I_2$.
\end{Lemma}
\begin{Proof}
The translation maps identify $E= \Jac(C)$ as a normal subgroup
of $\Aut(C)$. Conjugation by $g$ acts on $\Aut(C)$ and hence on $E$.
But the condition on the $j$-invariant ensures that the only
automorphisms of $E$ are $[\pm 1]$.
\end{Proof}

\begin{Lemma}
\label{normaliser2}
There is an exact sequence 
$$ \xymatrix{ 0 \ar[r] & \Theta_n' \ar[r] & N_{\G_n}(H_n) \ar[r]^-{\pi} &
\SL_2(\Z/n\Z) \ar[r] & 0 } $$
where $\Theta_n' \subset \G_n$ is generated 
by $H_n$ and the centre of $\G_n$.
\end{Lemma}
\begin{Proof}
If $n=2,3$ then $\G_n = \Gm \times \GL_n$ and the lemma already follows
by Lemma~\ref{normaliser1}.
If $n=4,5$ then $\G_n = \GL_m \times \GL_n$ where $m=2,5$.
Projection onto the second factor gives a 
map $\iota : N_{\G_n}(H_n) \to N_{\GL_n}(H_n)$
whose kernel is contained in the centre of $\G_n$.
The lemma will follow by Lemma~\ref{normaliser1} once we show that 
$\iota$ is surjective.

Let $B \in \GL_n$ with $B H_n B^{-1} = H_n$ and let $\phi \in X_n$
be a non-singular Hesse model. 
We know by Lemma~\ref{heiscond} that $C_\phi$ has Heisenberg group $H_n$.
If $\phi' = [I_m,B] \phi$ then 
$C_{\phi'}$ also has Heisenberg group $H_n$.
So by Lemmas~\ref{heis2} and~\ref{existsA}
there is a Hesse model $\phi''$ with $\phi'' = [A,I_n] \phi'$ 
for some $A \in \GL_m$. Putting $g=[A,B]$ we have
$\phi'' = g \phi$. Finally 
Lemma~\ref{heiscond} shows that $g H_n g^{-1} = H_n$.
\end{Proof}

Let $\alpha \in \SL_2(\Z/n\Z)$. By Lemma~\ref{normaliser2}
there exists $g \in N_{\G_n}(H_n)$ with 
$\pi(g) = \alpha$. Proposition~\ref{hn1}(ii) shows that $g$ acts on the
image of $u_n$. Since $u_n$ is linear we have 
$g \circ u_n = u_n \circ \gamma$ for some $\gamma \in \GL_2$.
Sending $\alpha$ to the class of $\gamma$ defines a group homomorphism
$$ \nu : \SL_2(\Z/n\Z) \to \PGL_2. $$
It is well defined by Lemmas~\ref{hcalc} and~\ref{normaliser2}.
The subgroup $\Delta_n \subset \PGL_2$ was 
defined in Section~\ref{sec:discinv}.

\begin{Lemma}
\label{propnu}
The map $\nu$ has kernel $\{ \pm I_2 \}$ and image $\Delta_n$.
\end{Lemma}

\begin{Proof}
It is possible to prove the lemma by a direct
calculation. An alternative method is as follows. 
First we use Lemma~\ref{pm1}
to show that the kernel of $\nu$ is contained in $\{ \pm I_2 \}$.
Then by Theorem~\ref{invthm}(ii) and Lemma~\ref{theymatch}
the image of $\nu$ permutes the roots of $D$.
Splitting into the cases $n=2,3,4,5$ it is easy to check
that $\Delta_n$ is the full group of such automorphisms
and $|\Delta_n| = |\PSL_2(\Z/n\Z)|$. The lemma follows
by counting.
\end{Proof}

\begin{ProofOf}{Proposition~\ref{hlem}}
(i) This can be proved by simply writing down a suitable element in 
the cases $n=2,3,4,5$.  An alternative method is as follows. 
Let $(a:b)$ be a point on $\PP^1$ that is fixed by no non-trivial
element of $\Delta_n$. Then $\phi = u_n(a,b)$ is a non-singular
Hesse model. We claim that there exists $g \in \G_n$ with 
$g \phi = \phi$ and $\det g = -1$. 
To prove this we first use \cite[Proposition 4.6]{g1inv}
to reduce to the case of a Weierstrass model, and then 
take $g= \gamma_n([-1;0,0,0])$ in \cite[Proposition 4.7]{g1inv}.
By Lemma~\ref{heiscond} we have $g H_n g^{-1} = H_n$
and so $g \circ u_n =u_n \circ \gamma$ for some $\gamma \in \GL_2$.
The image of $\gamma$ in $\PGL_2$ permutes the roots of $D$
and hence belongs to $\Delta_n$. Our choice of $(a:b)$ 
now forces $\gamma$ to be a scalar matrix. Since $g \phi = \phi$ it follows
that $g \circ u_n = u_n$. \\
(ii) We recall that $\GGamma_n$ is the inverse image of $\Delta_n$ in
$\SL_2$.  So this is immediate from Lemma~\ref{propnu}.
\end{ProofOf}

\begin{Remark}
It is possible to interpret $\nu : \SL_2(\Z/n\Z) \to \PGL_2$ 
as describing the automorphisms of the modular curve 
$X(n) \isom \PP^1$ obtained by relabelling the $n$-torsion
of the elliptic curves parametrised by $Y(n)$.
\end{Remark}

\begin{Remark}
By Definition~\ref{defGamma} and Lemma~\ref{propnu} 
both $\GGamma_n$ and $\SL_2(\Z/n\Z)$ are central extensions
of $\Delta_n$ by $\{ \pm 1 \}$.
In the cases $n=3,5$ we have $\GGamma_n \isom \SL_2(\Z/n\Z)$, but this is
not true for $n=2,4$.
\end{Remark}

\section{The Hesse polynomials}
\label{sec:hpolys}

The Hessian $H : X_n \to X_n$ was defined in Section~\ref{sec:covcontra}.

\begin{Lemma}
\label{stablem}
If $\phi \in X_n$ is non-singular then the subspace 
of $X_n$ fixed by 
the stabiliser of $\phi$ in $G_n$ 
is spanned by $\phi$ and $H(\phi)$.
\end{Lemma}
\begin{Proof}
By Proposition~\ref{hcovj} it suffices to prove the lemma 
for $\phi$ a Hesse model. Then by Lemma~\ref{heiscond} the stabiliser
is $H_n$ and by Proposition~\ref{hn1}(ii) the fixed subspace is
the Hesse family. Writing $\phi = u_n(a,b)$ it only
remains to check that $(a,b)$ and $(-\frac{\partial D}{\partial b},
\frac{\partial D}{\partial a})$ are linearly independent.
Since $D(a,b)^n = \Delta (\phi) \not= 0$ this is clear by
Euler's identity.
\end{Proof}

The pencil spanned by $U$ and $H$ has the following interpretation.

\begin{Theorem}
\label{hp}
Let $C \to \PP^{n-1}$ be a genus one normal curve
of degree $n=2,3,4,5$. Then $C=C_\phi$ for some $\phi \in X_n$ and 
\begin{enumerate}
\item If $\phi' = \lambda \phi + \mu H(\phi)$ is non-singular
then $C_{\phi'} \to \PP^{n-1}$ has the same Heisenberg group
as $C \to \PP^{n-1}$. 
\item If $C' \to \PP^{n-1}$ is a genus one normal
curve with the same Heisenberg group as $C \to \PP^{n-1}$ then
$C' = C_{\phi'}$ for some $\phi' = \lambda \phi + \mu H(\phi)$.
\end{enumerate}
\end{Theorem}

\begin{Proof} The existence of $\phi$ is clear for $n=2,3,4$.
In the case $n=5$ the genus one model $\phi$ is computed from the equations 
defining $C$ using the algorithm in \cite{g1pf}, based on the 
Buchsbaum-Eisenbud structure theorem \cite{BE1}, \cite{BE2} 
for Gorenstein ideals of 
codimension~$3$. \\
(i) By Proposition~\ref{hcovj} we may assume that $\phi$ 
is a Hesse model. Since the Hessian $H$
acts on the Hesse family it follows by Lemma~\ref{heiscond}
that both $C \to \PP^{n-1}$ and $C' \to \PP^{n-1}$ have 
Heisenberg group $H_n$ \\
(ii) Again we may assume that $\phi$ is a Hesse model.
Then Lemma~\ref{heis2} gives $C' = C_{\phi'}$ for 
some Hesse model $\phi'$. Since $\phi, \phi' \in X_n$ both have
stabiliser $H_n \subset G_n$ it follows by Lemma~\ref{stablem}
that $\phi' = \lambda \phi + \mu H(\phi)$ for some 
$\lambda, \mu \in \Kbar$.
\end{Proof}

In the later sections of this paper we will be concerned with
the arithmetic application of Theorem~\ref{hp}. First however
we record some formulae.

\begin{Lemma}
There are polynomials $f(\lambda,\mu)$, $g(\lambda,\mu)$
with coefficients in $K[c_4,c_6]$ such that
\[ H( \lambda U + \mu H) = f(\lambda,\mu)U + g(\lambda,\mu)H. \]
\end{Lemma}
\begin{Proof}
Writing $H( \lambda U + \mu H) = \sum F_{ij} \lambda^i \mu^j$
it is clear that the $F_{ij}$ are covariants. By 
Theorem~\ref{mainthm}(i) they are $K[c_4,c_6]$-linear combinations
of $U$ and $H$.
\end{Proof}

We compute the polynomials $f(\lambda,\mu)$ and $g(\lambda,\mu)$ 
just by working with the Hesse family. 
The case $n=3$ is classical: see  \cite[Section II.7]{Hilbert} or 
\cite[Art. 225]{Salmon}.
We will see in Theorem~\ref{hident} below 
that $f(\lambda,\mu)$ and $g(\lambda,\mu)$ are scalar multiples of
the partial derivatives of
\[ \D(\lambda,\mu)  =  
\left| \begin{matrix} \smallskip
\lambda & \mu \\ f(\lambda,\mu) & g(\lambda,\mu) 
\end{matrix} \right|. \]
\begin{Lemma}
\label{hpolysexist}
There are polynomials $\D(\lambda,\mu)$, $\cc_4(\lambda,\mu)$ 
and $\cc_6(\lambda,\mu)$ with coefficients in $K[c_4,c_6]$ 
such that
\begin{align*} \medskip
\D(\lambda,\mu) & =  D(\lambda a - \mu 
\tfrac{\partial D}{\partial b},\lambda b + 
\mu \tfrac{\partial D}{\partial a})/D(a,b) \\ \medskip
\cc_4(\lambda,\mu) & =  c_4(\lambda a - 
\mu \tfrac{\partial D}{\partial b},\lambda b + 
\mu \tfrac{\partial D}{\partial a}) \\
\cc_6(\lambda,\mu) & =  c_6(\lambda a - 
\mu \tfrac{\partial D}{\partial b},\lambda b + 
\mu \tfrac{\partial D}{\partial a}). 
\end{align*}
\end{Lemma}
\begin{Proof}
The coefficients of $\D(\lambda,\mu)$ belong to $K[a,b]$ 
since if $D(a,b)=0$ then $(a:b) = ( -\frac{\partial D}{\partial b} :
\frac{\partial D}{\partial a})$. The same is already clear for
$\cc_4(\lambda,\mu)$ and $\cc_6(\lambda,\mu)$. We must show that
the coefficients belong to $K[c_4,c_6]$.

Let $r= \deg c_4 = 4n/(6-n)$. Putting $\cc_4(\lambda,\mu) = \sum_{j=0}^r 
f_j \lambda^{r-j} \mu^j$ we find that $f_j \in K[a,b]$ has
degree $2n(2+j)/(6-n)$ and satisfies $f_j \circ \gamma = \chi^{2+j} (\gamma)
f_j$ for all $\gamma \in \GGamma_n$. So $f_j$ is a discrete invariant
satisfying the conditions of Theorem~\ref{proto}. It therefore belongs
to $K[c_4,c_6]$. The other cases are similar.
\end{Proof}

We call $\D(\lambda,\mu)$, $\cc_4(\lambda,\mu)$ and $\cc_6(\lambda,\mu)$ the 
{\em Hesse polynomials}. They are easily computed from the description 
in Lemma~\ref{hpolysexist}.
In the cases $n=2,3,4,5$ we find
\begin{align*}  \medskip
\D(\lambda,\mu) & =  \lambda^3-3 c_4 \lambda \mu^2-2 c_6 \mu^3 \\  \medskip
\D(\lambda,\mu) & =  \lambda^4-6 c_4 \lambda^2 \mu^2 - 8 c_6 
\lambda \mu^3-3 c_4^2 \mu^4 \\
\D(\lambda,\mu) & =  \lambda^6 - 15 c_4 \lambda^4 \mu^2 -40 c_6 
\lambda^3 \mu^3 - 45 c_4^2 \lambda^2 \mu^4  -  \, 24 c_4 c_6 \lambda \mu^5 
+ (27 c_4^3 - 32 c_6^2) \mu^6 \\ \medskip 
\D(\lambda,\mu) & =  \lambda^{12} - 66 c_4 
\lambda^{10} \mu^2 - 440 c_6 \lambda^9 \mu^3 
- 1485 c_4^2 \lambda^8 \mu^4 \\ & ~\hspace{2em} - 3168 
c_4 c_6 \lambda^7 \mu^5 + (5940 c_4^3 - 10560 c_6^2) \lambda^6 \mu^6
- 4752 c_4^2 c_6 \lambda^5 \mu^7  \\ 
& ~\hspace{2em} -  \, (66825 c_4^4 - 63360 c_4 c_6^2) 
\lambda^4 \mu^8
- (142560 c_4^3 c_6 - 140800 c_6^3) \lambda^3 \mu^9  \\
& ~\hspace{2em} -  \, (133650 c_4^5 -
133056 c_4^2 c_6^2) \lambda^2 \mu^{10} 
- (61560 c_4^4 c_6 - 61440 c_4 c_6^3) \lambda \mu^{11}  \\
& ~\hspace{2em} +  \, 
(91125 c_4^6 - 193536 c_4^3 c_6^2 + 102400 c_6^4) \mu^{12}.
\end{align*}
The polynomials $\D(\lambda,\mu)$ share with the $D(a,b)$
the property that their roots are arranged as the vertices of one of 
the Platonic solids.
By~(\ref{c4def}) and~(\ref{c6def}) we have
\[ \cc_4 (\lambda,\mu) =  
\tfrac{-1}{(\deg D)^2 ((\deg D)-1)^2} \left| \begin{matrix} \smallskip
\frac{\partial^2 \D}{\partial \lambda^2}(\lambda,\mu) & 
\frac{\partial^2 \D}{\partial \lambda \partial \mu}(\lambda,\mu) \\
\frac{\partial^2 \D}{\partial \lambda \partial \mu}(\lambda,\mu) &
\frac{\partial^2 \D}{\partial \mu^2}(\lambda,\mu) 
\end{matrix} \right| \]
and
\[ \cc_6(\lambda,\mu)  =  
\tfrac{1}{\deg D \deg c_4} \left| \begin{matrix} \smallskip
\frac{\partial \D}{\partial \lambda}(\lambda,\mu) & 
\frac{\partial \D}{\partial \mu}(\lambda,\mu) \\
\frac{\partial \cc_4}{\partial \lambda}(\lambda,\mu) &
\frac{\partial \cc_4}{\partial \mu}(\lambda,\mu) 
\end{matrix} \right|. \]
The Hesse polynomials are related by
\begin{equation}
\label{syz}
 \cc_4(\lambda,\mu)^3- \cc_6(\lambda,\mu)^2 
= (c_4^3-c_6^2) \, \D(\lambda,\mu)^n. 
\end{equation}

\begin{Theorem} 
\label{hident}
There are identities
\begin{align*}
c_4(\lambda U + \mu H) & =  \cc_4(\lambda,\mu) \\
c_6(\lambda U + \mu H) & =  \cc_6(\lambda,\mu) \\
(\deg D) H(\lambda U + \mu H) & =   
-\tfrac{\partial \D}{\partial \mu}(\lambda,\mu) 
U + \tfrac{\partial \D}{\partial \lambda}(\lambda,\mu) H. 
\end{align*}
\end{Theorem}

\begin{Proof}
As usual it suffices to check these relations on the Hesse family.
Each is a straightforward consequence of Lemma~\ref{hpolysexist}.
\end{Proof}

\begin{Remark}
We have $\Delta(\lambda U + \mu H) =  \Delta  \D(\lambda,\mu)^n$.
So the pencil spanned by 
a non-singular model and its Hessian has singular fibres at
the roots of $\D(\lambda,\mu)=0$. These may also be characterised
as the fibres whose Hessian is a scalar multiple of the original model.
(See also Theorem~\ref{thm:syz}.)
\end{Remark}

\section{The dual Hesse polynomials}
\label{sec:hpolys2}

In Section~\ref{sec:hpolys} we worked only with the invariants $c_4$ and $c_6$
and covariants $U$ and $H$. If we bring the contravariants
$P$ and $Q$ into play then there are many more identities to
consider. Again these are already in \cite{Salmon} in the case $n=3$.

\begin{Theorem} 
\label{dualthm1}
There are identities
\begin{align*}
 (\deg c_4) P(\lambda U + \mu H) & =  
f_4(\lambda,\mu)  P + g_4(\lambda,\mu) Q  \\
 (\deg c_6) Q(\lambda U + \mu H) & =   
f_6(\lambda,\mu) P + g_6(\lambda,\mu) Q  
\end{align*}
where
$$ \begin{pmatrix}
 f_4(\lambda,\mu) &  f_6(\lambda,\mu) \\
 g_4(\lambda,\mu) &  g_6(\lambda,\mu) 
\end{pmatrix}
= \begin{pmatrix}
 c_4 & c_6 \\ c_6 & c_4^2 
\end{pmatrix}^{-1}
\begin{pmatrix}
\frac{\partial \cc_4}{\partial \lambda}(\lambda,\mu) &  
\frac{\partial \cc_6}{\partial \lambda}(\lambda,\mu) \\
\frac{\partial \cc_4}{\partial \mu}(\lambda,\mu) &
\frac{\partial \cc_6}{\partial \mu}(\lambda,\mu) 
\end{pmatrix}. $$
\end{Theorem}

\begin{Proof}
Reducing to the Hesse family we must show that
$$ (\deg \Delta) \begin{pmatrix}
\partial c_4( \lambda U + \mu \partial D) \\
\partial c_6( \lambda U + \mu \partial D) 
\end{pmatrix}  =
 \begin{pmatrix}
 f_4(\lambda,\mu) &  g_4(\lambda,\mu) \\
 f_6(\lambda,\mu) &  g_6(\lambda,\mu) 
\end{pmatrix}
\begin{pmatrix} 3 \partial c_4 \\ 2 \partial c_6 \end{pmatrix}. $$
This follows from 
$$ 
(\deg D) D \begin{pmatrix}
\partial c_4( \lambda U + \mu \partial D) \\
\partial c_6( \lambda U + \mu \partial D) 
\end{pmatrix} = 
\begin{pmatrix}
\frac{\partial \cc_4}{\partial \lambda}(\lambda,\mu) &  
\frac{\partial \cc_4}{\partial \mu}(\lambda,\mu) \\
\frac{\partial \cc_6}{\partial \lambda}(\lambda,\mu) &
\frac{\partial \cc_6}{\partial \mu}(\lambda,\mu) 
\end{pmatrix}
\begin{pmatrix} \partial D \\ -U \end{pmatrix} $$
which is obtained by differentiating the definition in 
Lemma~\ref{hpolysexist}, and
$$ 1728 n D^{n-1} \begin{pmatrix} \partial D \\ -U \end{pmatrix}
= \begin{pmatrix}
 c_4^2 & -c_6 \\ -c_6 & c_4 
\end{pmatrix}
\begin{pmatrix} 3 \partial c_4 \\ 2 \partial c_6 \end{pmatrix} $$
which is a restatement of~(\ref{eq4}).
\end{Proof}

The {\em dual Hesse polynomials} $\DD(\xi,\eta)$, $\cdual_4(\xi,\eta)$, 
$\cdual_6(\xi,\eta)$ are defined in terms of the Hesse polynomials 
$\D(\lambda,\mu)$, $\cc_4(\lambda,\mu)$, $\cc_6(\lambda,\mu)$ by
\[ \begin{array}{lcl} \medskip
n = 2 & \quad & \begin{aligned}
\D(\lambda,\mu)  & =  -(c_4^3-c_6^2) \, \cdual_6(\xi,\eta) \\
\cc_4( \lambda,\mu) & =  (c_4^3-c_6^2)\, \cdual_4(\xi,\eta) \\
\cc_6( \lambda,\mu) & =  (c_4^3-c_6^2)^2 \, \DD(\xi,\eta) \\
\end{aligned} \\ \medskip
n = 3 & & \qquad \begin{aligned}
\D( \lambda,\mu)  & =  -(c_4^3-c_6^2) \, \cdual_4(\xi,\eta) \\
\cc_4( \lambda,\mu) & =  -(c_4^3-c_6^2)^2 \, \DD(\xi,\eta) \\
\cc_6( \lambda,\mu) & =  -(c_4^3-c_6^2)^2 \, \cdual_6(\xi,\eta) \\
\end{aligned} \\ \medskip
n = 4 & & \begin{aligned}
\D( \lambda,\mu)  & =  (c_4^3-c_6^2)^2 \, \DD(\xi,\eta) \\
\cc_4( \lambda,\mu) & =  (c_4^3-c_6^2)^2\,  \cdual_4(\xi,\eta) \\
\cc_6( \lambda,\mu) & =  (c_4^3-c_6^2)^3\,  \cdual_6(\xi,\eta) \\
\end{aligned} \\ \medskip
n = 5 & & \qquad \begin{aligned}
\D( \lambda,\mu)  & =  (c_4^3-c_6^2)^3 \, \DD(\xi,\eta) \\
\cc_4( \lambda,\mu) & =  (c_4^3-c_6^2)^4 \, \cdual_4(\xi,\eta) \\
\cc_6( \lambda,\mu) & =  (c_4^3-c_6^2)^6 \, \cdual_6(\xi,\eta) \\
\end{aligned}
\end{array} 
\]
where $\lambda = c_6 \xi + c_4^2 \eta$ and $\mu = -c_4 \xi - c_6 \eta$.
It follows by~(\ref{syz}) that
\[ \cdual_4(\xi,\eta)^3- \cdual_6(\xi,\eta)^2 
= (c_4^3-c_6^2)^{n-1} \DD(\xi,\eta)^n. \]

\begin{Theorem} 
\label{hidentdual}
There are identities
\begin{align*}
\tau^2 c_4(\xi P + \eta Q) & =  \cdual_4(\xi,\eta) \\
\tau^3 c_6(\xi P + \eta Q) & =  \cdual_6(\xi,\eta) \\
\tau (\deg D) H(\xi P + \eta Q) & =   
-\tfrac{\partial \DD}{\partial \eta}(\xi,\eta) P 
+ \tfrac{\partial \DD}{\partial \xi}(\xi,\eta) Q \\
\tau^2 (\deg c_4) P(\xi P + \eta Q) & =   
f_4(\xi,\eta)  U + g_4 (\xi,\eta)H  \\
\tau^3 (\deg c_6) Q(\xi P + \eta Q) & =   
f_6(\xi,\eta) U + g_6(\xi,\eta) H 
\end{align*}
where $\tau = 1,2,12,12^4$ for $n=2,3,4,5$ and
$$ \begin{pmatrix}
 f_4(\xi,\eta) &  f_6(\xi,\eta) \\
 g_4(\xi,\eta) &  g_6(\xi,\eta) 
\end{pmatrix}
= \begin{pmatrix}
 c_4 & c_6 \\ c_6 & c_4^2 
\end{pmatrix}^{-1}
\begin{pmatrix}
\frac{\partial \cdual_4}{\partial \xi}(\xi,\eta) &  
\frac{\partial \cdual_6}{\partial \xi}(\xi,\eta) \\
\frac{\partial \cdual_4}{\partial \eta}(\xi,\eta) &
\frac{\partial \cdual_6}{\partial \eta}(\xi,\eta)
\end{pmatrix}. $$
\end{Theorem}
\begin{Proof}
Again it suffices to check these identities on the Hesse family.
We did this by direct computation using {\sf MAGMA} \cite{magma}.
\end{Proof}
In the case $n=2$, the Hesse polynomials
and dual Hesse polynomials are the same. 
In {\sf MAGMA} our function
{\tt HessePolynomials(n,r,[c4,c6])} returns the 
Hesse polynomials $\D, \cc_4, \cc_6$ if $r=1$ and the dual 
Hesse polynomials $\DD, \cdual_4, \cdual_6$ if $r=-1$.

\section{Formulae}
\label{sec:formulae}

In the cases $n=2,3,4$ we give formulae for the
Hessian $H$ and for the contravariants $P$ and $Q$.
Theorem~\ref{kappa} then gives a practical method
for computing the invariants.
Alternatively we can compute the invariants using
the formulae in \cite{Mc+} or \cite[Section~7]{g1inv}.

\subsection{Formulae in the case $n=2$}

The binary quartic
$$f(x,z)  = a x^4 + b x^3 z + c x^2 z^2 + d x z^3 + e z^4 $$
has Hessian
\begin{align*} \medskip
H & =  (1/3) \times \left| \begin{matrix} \medskip
\frac{\partial^2 f}{\partial x^2} &  
\frac{\partial^2 f}{\partial x \partial z} \\
\frac{\partial^2 f}{\partial x \partial z} & 
\frac{\partial^2 f}{\partial z^2} 
\end{matrix} \right| \\
& =  (8 a c - 3 b^2) x^4 + (24 a d - 4 b c) x^3 z 
+ (48 a e + 6 b d - 4 c^2) x^2 z^2  \\ &  
\, ~\quad \,\,  + (24 b e - 4 c d) x z^3 + (8 c e - 3 d^2) z^4. 
\end{align*}
and contravariants
\begin{align*}  \medskip P & = 
e x^4 - d x^3 z + c x^2 z^2 - b x z^3 + a z^4 \\
Q & = (8 c e - 3 d^2) x^4 - (24 b e - 4 c d) x^3 z 
+ (48 a e + 6 b d - 4 c^2) x^2 z^2 \\ &  
\, \quad \,\, - (24 a d - 4 b c) x z^3 + (8 a c - 3 b^2) z^4.
\end{align*} 
The covariants and contravariants are closely related, 
the reason being that $X_2$ is isomorphic to its dual $X_2^*$ as
a $G_2= \SL_2$-module.

\subsection{Formulae in the case $n=3$}
The ternary cubic $U = U(x,y,z)$ has Hessian 
$$ \begin{array}{rcl} 
H & = & (-1/2) \times \left| \begin{matrix} \medskip
\frac{\partial^2 U}{\partial x^2} &  
\frac{\partial^2 U}{\partial x \partial y} &
\frac{\partial^2 U}{\partial x \partial z} \\ \medskip
\frac{\partial^2 U}{\partial x \partial y} & 
\frac{\partial^2 U}{\partial y^2} &
\frac{\partial^2 U}{\partial y \partial z} \\
\frac{\partial^2 U}{\partial x \partial z} &
\frac{\partial^2 U}{\partial y \partial z} & 
\frac{\partial^2 U}{\partial z^2} 
\end{matrix} \right|. 
\end{array} $$
The contravariant $P$, called in \cite{Salmon} the Caylean, is given by
$$ \begin{array}{rcl}
P  & = & (-1/xyz) \times \left| \begin{matrix} \medskip
\frac{\partial U}{\partial x} (0,z,-y) &  
\frac{\partial U}{\partial y} (0,z,-y) &
\frac{\partial U}{\partial z} (0,z,-y) \\ \medskip
\frac{\partial U}{\partial x} (-z,0,x) & 
\frac{\partial U}{\partial y} (-z,0,x) &
\frac{\partial U}{\partial z} (-z,0,x) \\
\frac{\partial U}{\partial x} (y,-x,0) &
\frac{\partial U}{\partial y} (y,-x,0) & 
\frac{\partial U}{\partial z} (y,-x,0) 
\end{matrix} \right|. 
\end{array} $$
The contravariant $Q$ may be computed from the 
coefficient of $\lambda^2 \mu$ in the first identity of
Theorem~\ref{dualthm1}, which in this case reads
$$ P( \lambda U + \mu H) = (\lambda^3 + 3 c_4 \lambda \mu^2 + 4 c_6 \mu^3) P 
 + 3( \lambda^2 \mu - c_4 \mu^3) Q. $$

\subsection{Formulae in the case $n=4$}

We identify a genus one model of degree $4$ with a pair of 
$4 \times 4$ symmetric matrices. Explicitly
$$ \phi = \begin{pmatrix} q_1 \\ q_2 \end{pmatrix} 
\equiv \begin{pmatrix} A \\ B \end{pmatrix} $$
where $q_1(x_1, \ldots, x_4) = \tfrac{1}{2} \x^T A \x$ 
and $q_2(x_1, \ldots, x_4)  =  \tfrac{1}{2} \x^T B \x.$
In the classical literature, as surveyed in \cite{Mc+}, the covariants
and contravariants are $\SL_4$-equivariant maps from 
$X_4$ to a space of quadrics.
In this setting the invariants $a, b, c, d, e$, 
contravariants $S_0, S_1, S_2, S_3$ and covariants
$A, T_1, T_2, B$ are given by
\begin{align*}
\det(sA+tB) & =  a s^4 + b s^3 t + c s^2 t^2 + d s t^3+e t^4 \\ 
\adj(sA+tB) & =  S_0 s^3+  S_1 s^2t + S_2 st^2 + S_3 t^3 \\
\adj(s(\adj A)+t(\adj B)) & =  a^2 A s^3 + a T_1 s^2 t + e T_2 s t^2 + e^2 B t^3.
\end{align*}
In terms of these, the Hessian is
\[ H = \begin{pmatrix}  6 T_2 - c A - 6 b B \\ 6 T_1 - c B - 6 d A
\end{pmatrix} \]
and the contravariants are
\[ P = \begin{pmatrix} \,\, 6e S_0 - 3 d S_1 + c S_2 -3 b S_3 \\
-3d S_0 + c S_1 -3 b S_2 + 6 a S_3 \end{pmatrix} \]
and 
\[  Q  =  \begin{pmatrix}
 (12 c e-18 d^2) S_0 + (- 18 b e+3 c d)S_1 + (12 a e+6 b d-c^2)S_2 
+ (-18 a d+ 3 b c) S_3 \\
 (-18 b e+3c d)S_0 + (12 a e+6 b d-c^2)S_1 + (-18 a d+3b c)S_2 + 
 (12 a c-18 b^2) S_3 \end{pmatrix}. \]

\section{An evaluation algorithm}
\label{sec:evaluation}

In the case $n=5$ the invariants $c_4$ and $c_6$ are homogeneous polynomials
of degrees $20$ and $30$ in $50$ variables. They are therefore
too large to compute as explicit polynomials. Nonetheless
we have found a practical algorithm for 
evaluating them (see \cite[Section~8]{g1inv}).
The Hessian $H : X_5 \to X_5$ is 
a $50$-tuple of homogeneous polynomials of degree $11$ 
in $50$ variables. Rather than attempt to compute these polynomials,
we shall again give an evaluation algorithm.

We identify $X_5 = \wedge^2 V \otimes W$ where $V$ and $W$ are 
$5$-dimensional vector spaces. Explicitly 
$$ ( \phi_{ij}(x_1, \ldots, x_5))_{i,j = 1, \ldots,5} \equiv
\sum_{i<j}  (v_i \wedge v_j) \otimes \phi_{ij}(x_1, \ldots, x_5) $$
where $v_1, \ldots, v_5$ and $x_1, \ldots,x_5$ are bases
for $V$ and $W$. The action of $\G_5 = \GL(V) \times \GL(W)$ is
the natural one. 
The covariants and contravariants considered so far are
the special cases $Y= X_5$ and $Y= X_5^*$ of the following
more general definition.

\begin{Definition}
Let $(\rho,Y)$ be a rational representation of $\G_5$. A covariant
is a polynomial map $F : \wedge^2 V \otimes W \to Y$ such that
$F \circ g= \rho(g) \circ F$ for all $g \in G_5.$
\end{Definition}

The $4 \times 4$ Pfaffians of $\phi \in X_5$ are quadrics $p_1, \ldots,p_5$ 
satisfying 
\[ \phi \wedge \phi \wedge v_i = p_i(x_1, \ldots, x_5) \,\, v_1 \wedge \ldots
\wedge v_5. \]
Let $v_1^*, \ldots, v_5^*$ be the basis for $V^*$ dual to 
$v_1, \ldots, v_5$. We define covariants
$$ \begin{array}{ll}
P_{2} : \wedge^2 V \otimes W \to V^* \otimes S^2W \, ; & \phi \mapsto 
\sum_{i=1}^5 v_i^* \otimes p_i(x_1, \ldots, x_5) \\
Q_{6} : \wedge^2 V \otimes W \to S^2 V \otimes W; & \phi \mapsto 
\sum_{i=1}^5 q_i(v_1, \ldots, v_5) \otimes x_i. \\
R_{10} :  \wedge^2 V \otimes W \to S^5 V^* ; & \phi \mapsto 
\det( \sum_{k=1}^5 
\frac{\partial^2 p_k}{\partial x_i \partial x_j} v_k^* ) \\
S_{10} : \wedge^2V \otimes W \to S^5 W \, ;  & 
  \phi \mapsto \det (\frac{\partial p_i}{\partial x_j})
\end{array} $$
where the auxiliary quadrics $q_i$ satisfy
\begin{equation}
\label{auxprop}
\tfrac{\partial}{\partial x_i} S_{10}(\phi)  =  q_i(p_1, \ldots, p_5).
\end{equation}
The proof that $Q_6$ exists, 
and is uniquely determined by~(\ref{auxprop}), 
is given in \cite[Section 8]{g1inv}.
We write $\langle~,~\rangle$ for the contraction
\begin{align*}
S^a V \times S^{a+b} V^* & \to  S^{b} V^*  \\
(f(v_1, \ldots, v_5),g(v_1^*, \ldots, v_5^*))
& \mapsto  f(\tfrac{\partial }{\partial v_1^*}, \ldots,
\tfrac{\partial }{\partial v_5^*}) g(v_1^*, \ldots, v_5^*) 
\end{align*}
and identify $X_5 = \wedge^2 V \otimes W$ with the 
space of $10 \times 5$ matrices via 
$$ \sum_{i<j} v_i \wedge v_j \sum_k a_{ijk} x_k \equiv \begin{pmatrix} 
a_{121} & a_{122} & \cdots &  a_{125} \\
a_{131} & a_{132} & \cdots &  a_{135} \\
\vdots & \vdots  & & \vdots \\
a_{451} & a_{452} & \cdots &  a_{455} \\
 \end{pmatrix}. $$
\begin{Theorem} 
\label{hessalgthm}
The Hessian $H: X_5 \to X_5$ satisfies
\begin{align*}
P_2 \circ H  & =  
4 c_4 P_2 - \tfrac{3}{16} \langle Q_6 ,\langle Q_6, R_{10} \rangle 
\rangle \\ \det (U;H) & =  12^5 \Delta. 
\end{align*}
These conditions uniquely determine $H(\phi)$ for $\phi \in X_5$ 
non-singular.
\end{Theorem}
\begin{Proof}
The covariance of these identities is clear, so it suffices to check
them for $\phi$ a Hesse model.
We did this by direct calculation. If $\phi \in X_5$ is non-singular
then $H(\phi)$ is equivalent to a Hesse model, and therefore defines
a curve. So for the final statement all we need to know is that
if $\phi_1, \phi_2 \in X_5$ each define a curve
and $P_2 (\phi_1)= P_2(\phi_2)$ then $\phi_1= \pm \phi_2$.
This follows from the Buchsbaum-Eisenbud acyclicity criterion and
the properties of minimal free resolutions. See for example
\cite[Section 5.2]{g1inv}.
\end{Proof}

To compute the Hessian of a non-singular model 
$\phi \in X_5$ we begin by computing its invariants
using the algorithm in~\cite[Section~8]{g1inv}. 
The auxiliary quadrics $q_1, \ldots, q_5$ are
computed as a by-product of this algorithm. Then we use the 
first identity of Theorem~\ref{hessalgthm} to compute the $4 \times 4$
Pfaffians of $H(\phi)$. The genus one 
model $H(\phi)$ is recovered from its $4 \times 4$ Pfaffians 
using the algorithm in \cite{g1pf}. This only determines 
$H(\phi)$ up to sign. In applications where the sign 
matters we use the second identity in Theorem~\ref{hessalgthm} to make 
a consistent choice.

\begin{Remark}
\label{othercov}
We have found similar algorithms for computing the contravariants
(i.e. covariants for $Y = \wedge^2 V^* \otimes W^*$)
and also the covariants for $Y=\wedge^2 W \otimes V^*$
and $Y= \wedge^2 W^* \otimes V$. 
We will report on these constructions and their arithmetic applications
in \cite{enqI}. (See also Remark~\ref{other5}.)
\end{Remark}

\section{Theta groups}
\label{Theta groups}

Let $C \to \PP^{n-1}$ be a genus one normal curve defined over $K$.
We recall that $E =\Jac(C)$ acts on $C$ by translation, and 
translation by $P \in E$ extends to an automorphism of $\PP^{n-1}$
if and only if $P \in E[n]$. 
Analogous to Definition~\ref{heisdef} we have

\begin{Definition} 
\label{def:theta}
The {\em theta group} of $C \to \PP^{n-1}$ is the 
group of all matrices in $\GL_n$ that act on $C$ as translation
by some $T \in E[n]$. As a group it is a central extension of 
$E[n]$ by $\Gm$ with commutator given by the Weil pairing.
\end{Definition}

Let $T \mapsto M_T$ be a Galois equivariant section for $\Theta \to E[n]$.
We consider the following pair of inverse problems.

\begin{enumerate}
\item \label{problemi} 
Given equations for $C \to \PP^{n-1}$ how can we compute the matrices
$M_T$ for $T \in E[n]$? 
\item \label{problemii}
Given the matrices $M_T$ for $T \in E[n]$ how can we compute 
equations for $C \to \PP^{n-1}$?
\end{enumerate}

We restrict to $n =2,3,4,5$ and write $C = C_\phi$ where $\phi$ is 
a genus one model of degree $n$ defined over $K$.  
The pencil of curves spanned by $\phi$ 
and its Hessian is a twist of the Hesse family.
In particular there are $12/(6-n)$ singular fibres, and each of
these is an $n$-gon. Generalising the terminology in 
\cite[Section II.7]{Hilbert} 
we call these the
syzygetic $n$-gons. 
If we change co-ordinates so that one of the syzygetic $n$-gons
has vertices $(1:0:0: \ldots), \, (0:1:0: \ldots), \ldots$ then 
the action of some $T \in E[n]$ is given by $M_T = \Diag(1,\zeta_n,
\zeta_n^2, \ldots)$ where $\zeta_n$ is a primitive $n$th root of unity.
Conversely the following theorem gives formulae for a syzygetic
$n$-gon in terms of $\zeta_n$ and $T$. 

\begin{Theorem}
\label{thm:syz}
Let $\phi$ be a non-singular genus one model of degree $n = 2,3,4,5$
with invariants $c_4$ and $c_6$. Let $T = (x_T,y_T)$ be
a torsion point of order~$n$ on the Jacobian $E:  \,\, 
y^2 = x^3 - 27 c_4 x - 54 c_6$. Then there is 
a syzygetic $n$-gon defined by $\psi = \xi_T \phi + 3 H(\phi)$
and if $n \ge 3$ 
this model satisfies $H(\psi) = - (\frac{1}{3} \eta_T)^2 \psi$ where
\[  \xi_T = \left\{ \begin{array}{ll} x_T & \text{ if $n=2, 3, 4$ }  \\
(1+\zeta_5 + \zeta_5^4) x_T + (1+\zeta_5^2 + \zeta_5^3) x_{2T} &
\text{ if } n = 5 \end{array} \right. \]
and
\[  \eta_T = \left\{ \begin{array}{ll} (\zeta_3 - \zeta_3^{-1}) y_T 
& \text{ if } n = 3  \\
(\zeta_4 - \zeta_4^{-1}) (x_T - x_{2T}) y_T 
& \text{ if } n = 4  \\
(x_T - x_{2T})^4 \big((\zeta_5 - \zeta_5^4)^5 y_T 
+ (\zeta_5^2 - \zeta_5^3)^5 y_{2T} \big)
& \text{ if } n = 5. \end{array}  \right. \]
\end{Theorem}

\begin{Proof}
The syzygetic $n$-gons are the fibres of the pencil $\lambda \phi 
+ \mu H(\phi)$
where $(\lambda : \mu)$ is a root of Hesse polynomial $\D(\lambda,\mu)=0$.
A calculation using division polynomials shows that $\D(\xi_T,3)=0$.
Then by Theorem~\ref{hident} and a further calculation using 
division polynomials $H(\psi) = \tfrac{1}{3 (\deg \D)} \,
 \tfrac{\partial \D}{\partial \lambda} (\xi_T,3) \psi
  = - (\tfrac{1}{3} \eta_T)^2 \psi$.
\end{Proof}

\begin{Remark}
The syzygetic $n$-gon in Theorem~\ref{thm:syz} has field of
definition $K(\xi_T)$, whereas an orientation of the $n$-gon is defined
over $K(\xi_T,\eta_T)$. If $\zeta_n \in K$ then these fields are just
$K(x_T)$ and $K(T) = K(x_T,y_T)$. 
\end{Remark}

Theorem~\ref{thm:syz} is the basis for our method for computing the $M_T$,
i.e. for solving problem~(\ref{problemi}). We have worked out
explicit formulae in the cases $n=2,3,4$. These are given in 
\cite[Section 6.2]{minred234}, \cite{testeqtc} and~\cite{improve4}. 
Applications include testing equivalence of genus one 
models, adding Selmer group elements (represented as explicit
covering curves), and computing the inner product used 
for the reduction of genus one models. 
In the case $n=5$ we have not yet found formulae for the $M_T$ 
as explicit as the ones we gave for $n=2,3,4$.
 
We now turn to problem~(\ref{problemii}), i.e. we try to recover $C$
from the $M_T$. This is required for the Hesse pencil method of 
$n$-descent as described in \cite[Section 5.1]{paperI}. 
See also the example of $5$-descent in \cite[Section 9.3]{paperIII}. 

We return to considering arbitrary $n \ge 2$. Generalising 
Definition~\ref{def:theta} we have

\begin{Definition} 
\label{def:theta2}
A theta group for $E[n]$ is a central
extension of $E[n]$ by $\Gm$ with commutator given by the 
Weil pairing.
\end{Definition}

The maps $\Gm \to \Theta$ and $\Theta \to E[n]$
are considered part of the data defining the theta group. Thus an 
isomorphism of theta groups is a commutative diagram
\[ \xymatrix{ 0 \ar[r] & \Gm \ar[r] \ar@{=}[d] & \Theta_1 \ar[r] 
\ar[d]^{\isom} & E[n] \ar[r] \ar@{=}[d] & 0 \\ 0 \ar[r] & \Gm \ar[r]
 & \Theta_2 \ar[r] & E[n] \ar[r] & 0  } \]
If we ignore the Galois action (i.e. by passing to the algebraic closure
$\Kbar)$, then there is only one such group. Its automorphism group 
is naturally a copy of $\Hom(E[n],\Gm)$ which we identify with $E[n]$ 
via the Weil pairing. Let $\Theta_E$ 
be the base choice of theta group arising from the construction of
Definition~\ref{def:theta} for $E$ embedded in $\PP^{n-1}$ via the 
complete linear system $|n.0_E|$. Then as noted in \cite[Section 1.6]{paperI} 
the group $H^1(K,E[n])$ parametrises the theta groups for $E[n]$ (up to 
$K$-isomorphism) as twists of $\Theta_E$.

\begin{Theorem}
\label{xn-family}
Let $E[n] \to \GL_n\,; \,\, T \mapsto M_T$ 
be a Galois equivariant map arising as a section of a theta 
group $\Theta$ for $E[n]$. Then 
\begin{enumerate}
\item There is a genus one normal curve $C \to \PP^{n-1}$ 
defined over $K$ with Jacobian $E$ and theta group $\Theta$. 
\item If $n \ge 3$ then the genus one normal 
curves  $C \subset \PP^{n-1}$ for which
each matrix $M_T$ acts as translation by some $n$-torsion point
of $\Jac(C)$ are parametrised by a twist of the modular
curve $Y(n)$. 
\end{enumerate}
\end{Theorem}
\begin{Proof}
See \cite[Theorem 5.2 and Proposition 5.5]{paperI}. 
\end{Proof}

Theorem~\ref{xn-family}(ii) is false in the case $n=2$
since replacing a $2$-covering $y^2 = g(x)$
by a quadratic twist $d y^2 = g(x)$ does not change
the matrices $M_T$.

Again restricting to $n = 2,3,4,5$ we explain how to compute the
family of curves in Theorem~\ref{xn-family}(ii). In the cases $n=2,3$
we solve for the $2$-dimensional Heisenberg invariant subspace
of $X_n$. The case $n=3$ is also described in \cite[Lemma 5.6]{paperI}. 
The cases $n=4$, $5$ are
more complicated since we do not yet know the action of $E[n]$ on the
space of equations defining $C \subset \PP^{n-1}$. 

\begin{Lemma} 
\label{lem:4quads}
If $n=4$ then the space of quadrics $Q \in K[x_1, \ldots, x_4]$ satisfying
$Q \circ M_T^2 = - (\det M_T) Q$ for all $T \in E[4] \setminus E[2]$ 
has dimension $4$.
\end{Lemma}
\begin{Proof} We first note that the condition on the quadrics 
is independent of the scaling of the matrices $M_T$. 
Working over $\Kbar$ we may assume 
that $E[4]$ is generated by $T_1$ and $T_2$ acting via
\[ M_{T_1} = \left( \begin{smallmatrix} 
1 & & & \\ & i & & \\ & & -1 & \\ & & & -i
\end{smallmatrix} \right) \quad \text{ and } \quad  
M_{T_2} =  \left( \begin{smallmatrix} & & & 1 \\
1 & & & \\ & 1 & & \\ & & 1 & \end{smallmatrix} \right).  \]
Then the space of quadrics in question has basis
$x_1^2 + x_3^2$, $x_1 x_3$, $x_2^2 + x_4^2$, $x_2 x_4$.
\end{Proof}

Let $H \isom H_4$ be the intersection of $\Theta$ with $\SL_4$.
Lemma~\ref{lem:4quads} constructs a $4$-dimensional representation of 
$H$, say $V$. Intersecting the $H$-invariant subspace of $\wedge^2 V$
with the image of $V \times V \to \wedge^2 V ; \, (v_1,v_2) 
\mapsto v_1 \wedge v_2$ gives a conic $\Gamma \subset \PP(\wedge^2 V)$.
Each point on $\Gamma$ corresponds to a $2$-dimensional $H$-invariant
space of quadrics. We have thus constructed the family of curves
specified in Theorem~\ref{xn-family}(ii). Theorem~\ref{xn-family}(i) 
shows that $\Gamma(K) \not= \emptyset$ and hence $\Gamma \isom \PP^1$.

In the case $n=5$ we again let $H \isom H_5$ be the intersection of
$\Theta$ with $\SL_5$. The space of linear forms on $\PP^4$
is a $5$-dimensional representation of $H$, say $W$. If we ignore the
Galois action then this is the unique irreducible representation of $H$
with central character $\zeta \mapsto \zeta$. Again by inspection of 
the character table for $H_5$ there is a unique irreducible representation
$V$ of $H$ with central character $\zeta \mapsto \zeta^2$. 
Since $V$ is $5$-dimensional the representations
$\wedge^2 W$ and $S^2 W$ are isomorphic (over $\Kbar$) to either $2$
or $3$ copies of $V$. We solve by linear algebra for a non-zero
$H$-equivariant $K$-linear map $\pi : S^2 W \to \wedge^2 W$. Then 
either the kernel or image of $\pi$ is a copy of $V$. The proof of
Proposition~\ref{hn1} shows that the space $(\wedge^2 V \otimes W)^H$ is 
$2$-dimensional. Since our construction of $V$ from $W$ is Galois
equivariant, we can find a basis for this space defined over $K$.
Identifying $\wedge^2 V \otimes W$ with the space of genus one
models of degree $5$ we have thus constructed the family of curves
specified in Theorem~\ref{xn-family}(ii).

We return to considering the cases $n = 2,3,4,5$. Let $\phi$ be a non-singular
genus one model defining a curve in the family specified in 
Theorem~\ref{xn-family}(ii). By Theorem~\ref{hp} we can  
recover the whole family by taking linear combinations of $\phi$
and its Hessian. Problem~(\ref{problemii}), stated at the start of
this section, is now reduced~to

\begin{enumerate}
\setcounter{enumi}{2}
\item \label{problemiii} Given a non-singular genus one model
$\phi \in X_n(K)$, find all 
models $\phi' = \lambda  \phi + \mu H(\phi)$
with $(\lambda: \mu) \in \PP^1(K)$ and $\Jac(C_{\phi'}) \isom E$.
\end{enumerate}

Let $c_4,c_6,\Delta$ be the invariants of $\phi$ and $\D, \cc_4, \cc_6$ the
Hesse polynomials, with coefficients evaluated at $c_4, c_6 \in K$.
Then the fibres with the same $j$-invariant as $E$ correspond
to the roots $(\lambda:\mu) \in \PP^1$ of the binary form
\[ \cc_4(\lambda,\mu)^3 - j(E)  \, \Delta \, \D(\lambda,\mu)^n . \]
Solving for the $K$-rational roots of a polynomial of 
degree $12n/(6-n)$ we are left with only finitely many possibilities
for $\phi'$ (up to quadratic twists in the case $n=2$). 
We then use Theorem~\ref{invthm}(iii) to decide which 
of these define a curve with Jacobian $E$.

As noted in \cite[Section 5.1]{paperI} it can happen that the family
of curves in Theorem~\ref{xn-family}(ii) contains more than one
curve defined over $K$ with Jacobian~$E$. However this does not happen
in the generic case, i.e. when \[\Gal(\Kbar/K) \to \Aut(E[n]) \isom 
\GL_2(\Z/n\Z)\] is surjective.
In the event that more than one curve remains, we use our solution to 
problem~(\ref{problemi}) to determine which of these is correct.

\section{The modular curves $X_E(n)$ and $X^-_E(n)$}
\label{famec}


\begin{Definition}
\label{def:ncongr}
Elliptic curves $E$ and $E'$ defined over $K$ are  
{\em directly $n$-congruent} if there is an isomorphism of
Galois modules $E[n] \isom E'[n]$ that respects the Weil pairing.
They are {\em reverse $n$-congruent} if there 
is an isomorphism of Galois modules $\psi : E[n] \isom E'[n]$ satisfying 
$e_n(\psi S,\psi T) = e_n(S,T)^{-1}$ for all $S,T \in E[n]$.
\end{Definition}

The elliptic curves directly $n$-congruent to $E$ are parametrised
by $Y_E(n)$ and those reverse $n$-congruent to $E$ by $Y_E^-(n)$.
The smooth projective models $X_E(n)$ and $X_E^-(n)$ are twists of 
$X(n)$.


In the cases $n=2,3,4,5$ we have $X(n) \isom \PP^1$. We show that
the Hesse polynomials define the families
of curves parametrised by $X_E(n) \isom \PP^1$. Recall that in the 
case $n=3$ the Hesse polynomials were already known to Salmon 
\cite{Salmon}.


\begin{Theorem}
\label{firsttwist}
Let $n=2,3,4,5$.
Let $E$ be an elliptic curve over $K$,
$$ y^2 = x^3 - 27 c_4 x - 54 c_6,$$
and let $E_{\lambda,\mu}$ be the family of curves
$$ y^2 = x^3 - 27 \cc_4(\lambda,\mu) x - 54 \cc_6(\lambda,\mu)$$
where the coefficients of the Hesse polynomials 
$\cc_4(\lambda,\mu)$ and $\cc_6(\lambda,\mu)$ 
are evaluated at $c_4,c_6 \in K$. 
Then an elliptic curve $E'$ over $K$ is directly $n$-congruent to $E$
if and only if it is isomorphic over $K$ to $E_{\lambda,\mu}$
for some $\lambda, \mu \in K$.
\end{Theorem}
\begin{Proof}
Let $E \to \PP^{n-1}$ be the genus one normal curve given 
by the complete linear system $|n.0_E|$. It is defined by some 
$\phi \in X_n(K)$ with invariants $c_4$ and $c_6$. \\
(i) Suppose that $\phi' = \lambda \phi + \mu H(\phi)$ is non-singular. 
By Theorem~\ref{hp}(i)
the genus one normal curves $C_\phi \to \PP^{n-1}$ 
and $C_{\phi'} \to \PP^{n-1}$ have the same 
Heisenberg group. By Theorems~\ref{invthm}(iii) and~\ref{hident}
their Jacobians are $E$ and $E_{\lambda,\mu}$. It follows by
Definition~\ref{heisdef} that $E$ and $E_{\lambda,\mu}$ 
are $n$-congruent. \\
(ii) By Theorem~\ref{xn-family}(i) there is a genus one normal curve 
$C' \to \PP^{n-1}$ with Jacobian $E'$ and the same Heisenberg
group as $E \to \PP^{n-1}$. Then Theorem~\ref{hp}(ii)
shows that $C' = C_{\phi'}$ for some $\phi' = \lambda \phi + \mu H(\phi)$.
Since $C'$ is defined over $K$ we may arrange that 
$\lambda, \mu \in K$.
Taking Jacobians 
gives $E' \isom E_{\lambda,\mu}$.
\end{Proof}

If we split into the cases $c_4 c_6 \not=0$, $c_4=0$, $c_6=0$,
then Theorem~\ref{firsttwist} reduces to formulae 
obtained by Rubin and Silverberg \cite{RubinSilverberg},
\cite{RubinSilverberg2}, \cite{Silverberg}.
To explain the relationship in the case $c_4 c_6 \not=0$ we write
$$ H_1= \frac{-c_6^2 U + c_4 c_6 H}{c_4^3-c_6^2}. $$
\begin{Lemma}
\label{Jpolys}
There are polynomials $\alpha(J,t)$ and $\beta(J,t)$ such that
\begin{align*}
c_4(U + t H_1 ) & =  \cc_4( 1- \tfrac{c_6^2}{c_4^3-c_6^2}t,
\tfrac{c_4 c_6}{c_4^3-c_6^2}t )
  =  \alpha(J,t) c_4 \\
c_6(U + t H_1 ) & =  \cc_6( 1- \tfrac{c_6^2}{c_4^3-c_6^2}t,
\tfrac{c_4 c_6}{c_4^3-c_6^2}t )
  =  \beta(J,t) c_6 
\end{align*}
where $J= c_4^3/(c_4^3-c_6^2)$.
\end{Lemma}
\begin{Proof}
Theorem~\ref{hident} gives the first two equalities. Then expanding 
 $c_4(U+tH_1)/c_4$ and $c_6(U+tH_1)/c_6$ in powers of $t$ we see that
the coefficients are weight zero elements of $K[c_4,c_6,\Delta^{-1}]$
and therefore polynomials in $J$.
\end{Proof}

\begin{Corollary}
\label{rscor}
Let $n=3,4,5$.
Let $E$ be an elliptic curve over $K$,
$$ y^2 = x^3 + a x + b,$$
and let $E_t$ be the family of curves
$$ y^2 = x^3 + \alpha(J,t) a x + \beta(J,t) b $$
where $J= j(E)/1728 = 4a^3/(4a^3+27 b^2).$ \\
(i) Every elliptic curve $E_t$ over $K$ with $t \in \PP^1(K)$ 
is directly $n$-congruent to $E$. \\
(ii) If $j(E) \not= 0,1728$ then every elliptic curve $E'$ over $K$, 
that is directly $n$-congruent to $E$, 
is isomorphic over $K$ to $E_t$ for some $t \in \PP^1(K)$.
\end{Corollary}
\begin{Proof}
This follows from Theorem~\ref{firsttwist} and Lemma~\ref{Jpolys}. 
We assume $n=3,4,5$ so that $E_{\lambda,\mu}$ is determined
up to isomorphism by the ratio $(\lambda: \mu) \in \PP^1(K)$.
(This is false for $n=2$.)
The condition $j(E) \not =0,1728$ is required so that the matrix 
$$ \begin{pmatrix} 1 & \frac{-c_6^2}{c_4^3-c_6^2} \\
0 & \frac{c_4 c_6}{c_4^3-c_6^2} \end{pmatrix}
$$ is non-singular.
\end{Proof}

The polynomials $\alpha(J,t)$ and $\beta(J,t)$ 
for $n=3,5,4,2$ are written out in 
\cite[Theorem 4.1 and Appendix]{RubinSilverberg}, 
\cite{Silverberg}, \cite{RubinSilverberg2}.
They are also returned by our {\sf MAGMA} function
{\tt RubinSilverbergPolynomials}.

Our  proof of Theorem~\ref{firsttwist} only requires
the existence of the Hessian, not the explicit formulae and algorithms
in Sections~\ref{sec:formulae} and~\ref{sec:evaluation}. Alternatively
the theorem can be proved (without using the Hessian) by adapting 
the method of Rubin and Silverberg. For this we replace the
$2$ by $2$ matrix $A$ used in \cite{RubinSilverberg}
(which is singular for $j(E) = 0,1728$ --
hence the restriction in Corollary~\ref{rscor}(ii)) by the matrix
\[ \begin{pmatrix} a & -\frac{\partial D}{\partial b} \\ b  
& \,\,\, \frac{\partial D}{\partial a} \end{pmatrix} \]
used in Lemma~\ref{hpolysexist}.

The analogue of Theorem~\ref{firsttwist} for 
$X_E^-(n)$ 
is obtained by replacing the Hesse polynomials
$\cc_4(\lambda,\mu)$ and $\cc_6(\lambda,\mu)$ by the dual Hesse polynomials
$\tau^{-2} \cdual_4(\xi,\eta)$ and $\tau^{-3} \cdual_6(\xi,\eta)$.
The proof uses the contravariants $P$ and $Q$ instead of
the covariants $U$ and $H$. To see why this works, recall that if 
$M_T$ is a matrix describing the action of $T \in E[n]$ on $C \to \PP^{n-1}$
then the Weil pairing is given by $e_n(S,T)I_n = M_S M_T M_S^{-1} M_T^{-1}$.
Replacing each matrix $M_T$ by its inverse transpose therefore has the
effect of switching the sign of the Weil pairing.

If $n=2$ or $5$ then $n$-congruence and
reverse $n$-congruence are the same (since $-1$ is a square mod $n$).
The families parametrised by $X_E(n)$ and $X_E^-(n)$ are therefore
isomorphic. If $n=4$ then (by the formulae in Section~\ref{sec:hpolys2}) 
the only change is that we take the
quadratic twist by the discriminant $\Delta$. The analogue of 
Corollary~\ref{rscor}
for reverse $3$-congruence holds for the family of curves
\[ y^2 = x^3 - 4 \gamma(J,t) J a x - 8 \beta(J,t) J^2 b \]
where $\gamma(J,t) =  \D(1-\tfrac{c_6^2}{c_4^3-c_6^2}t,
\tfrac{c_4 c_6}{c_4^3-c_6^2} t).$ 
By~(\ref{syz}) and Lemma~\ref{Jpolys} 
we may characterise $\gamma(J,t)$ as the unique polynomial 
in $\Z[J,t]$ satisfying
$\gamma(J,t)^3 = \alpha(J,t)^3 J + \beta(J,t)^2 (1-J)$.

\section{Visibility of Tate-Shafarevich groups}
\label{sec:vis}

In this section we recall the theory of visibility, introduced by
Mazur \cite{CM}, \cite{Mazurcubics}, and explain how we can use
it to compute explicit elements of the Tate-Shafarevich group of 
an elliptic curve. We give some examples in Section~\ref{sec:examples}.

We start with a short exact sequence of abelian varieties
\[ 0 \ra E \stackrel{\iota}{\ra} A \stackrel{\psi}{\ra} F' \ra 0 \]
defined over a number field $K$. If $P \in F'(K)$ then restricting
the group law on $A$ gives the fibre $\psi^{-1}(P)$ the structure
of torsor under $E$. It therefore represents an element of the
Weil-Ch\^atelet group $H^1(K,E)$. In fact this element is $\delta(P)$
where $\delta$ is the connecting map in the long exact sequence
\[ \ldots \ra F'(K) \stackrel{\delta}{\ra} H^1(K,E) 
\stackrel{\iota_*}{\ra} H^1(K,A) \ra \ldots \]
Following Mazur we define $\xi \in H^1(K,E)$ to be {\em visible} in $A$ 
if $\xi = \delta(P)$ for some $P \in F'(K)$, equivalently 
$\iota_*(\xi) = 0$. 

Now let $E$ and $F$ be elliptic curves with a common 
finite Galois submodule $\Phi$. We put $E' = E/\Phi$ and $F' = F/\Phi$.
Then $A = (E \times F)/\Phi$ is an abelian surface and taking 
Galois cohomology of 
\[ \xymatrix{ 0 \ar[r] & \Phi \ar[r] \ar[d] & F \ar[d] \ar[r] &
F' \ar[r] \ar@{=}[d] & 0 \\ 0 \ar[r] & E \ar[r] & A \ar[r] 
& F' \ar[r] & 0  } \]
gives the following commutative diagram whose rows are the 
Kummer exact sequences for $\phi_E: E \to E'$ and $\phi_F : F \to F'$.
\begin{equation*}
\xymatrix{ 0 \ar[r] &  \frac{E'(K)}{\phi_E E(K)} 
\ar[r] & H^1(K,\Phi) \ar@{=}[d] \ar[r] &
H^1(K,E) [ \phi_E ]  \ar[r] & 0 \\ 0 \ar[r] &  
\frac{F'(K)}{\phi_F F(K)}\ar[r]  \ar[urr]_(0.6){\delta}|!{[r];[ur]}\hole 
& H^1(K,\Phi) \ar[r] 
& H^1(K,F)[\phi_F] \ar[r] & 0  } 
\end{equation*}

We are interested in the case $\Phi = E[n]= F[n]$ for some integer
$n \ge 2$. In particular $\phi_E$ and $\phi_F$ are multiplication-by-$n$
on $E=E'$ and $F=F'$. 

In \cite{paperI} we gave a list of geometric interpretations 
of the group $H^1(K,E[n])$. Two of these are relevant here. The first is that
$H^1(K,E[n])$ parametrises the torsor divisor class pairs $(C,[D])$
as twists of $(E,[n.0_E])$. The second, already recalled in 
Section~\ref{Theta groups}, is that $H^1(K,E[n])$ parametrises the 
theta groups for $E[n]$ as twists of $\Theta_E$. If $D$ is a $K$-rational
divisor then $(C,[D])$ determines a morphism 
$C \to \PP^{n-1}$ via the complete linear system $|D|$ and then a 
theta group by the construction of Definition~\ref{def:theta}. 
We checked in \cite[Section 1.6]{paperI} that this 
construction is compatible with the above two interpretations of
$H^1(K,E[n])$.

\begin{Lemma}
\label{lem:n_odd}
If $n$ is odd then $\Theta_E$ depends only on $E[n]$ (regarded as a
Galois module equipped with the Weil pairing) and not on $E$.
\end{Lemma}
\begin{Proof}
See \cite[Lemma 3.11]{paperI}. 
\end{Proof}

Our method for computing equations for visible elements of $\Sha(E/K)[n]$
is as follows. First we find a second elliptic curve $F/K$ that is
directly $n$-congruent to $E$ (see Definition~\ref{def:ncongr}).
Then we
convert $P \in F(K)$ to a torsor $C$ under $E$ by realising the map
$\delta$ in the above diagram as the composite
\begin{equation}
\label{eq:composite}
 \frac{F(K)}{n F(K)} \stackrel{\delta_F}{\ra} H^1(K,F[n]) 
\isom H^1(K,E[n]) \stackrel{\iota_{E,*}}{\ra} H^1(K,E)[n]. 
\end{equation}
We recall from \cite[Remark 1.10]{paperI} that, in terms of torsor divisor 
class pairs, the maps $\delta_F$ and $\iota_{E,*}$
are given by $P \mapsto (F,[(n-1).0_F + P])$ and $(C,[D]) \mapsto C$. 
In contrast we realise the middle isomorphism in~(\ref{eq:composite}) 
by using that $E[n]$ and $F[n]$ have the same theta groups 
(see Definition~\ref{def:theta2}).

We make this construction explicit in the cases $n = 2,3,4,5$. First 
we map $F \to \PP^{n-1}$ via the complete linear system $|(n-1).0_F + P|$,
and let $\phi \in X_n(K)$ be a genus one model defining the image.
Then we compute the Hessian $H(\phi)$ using the formulae and algorithms 
in Sections~\ref{sec:formulae} and~\ref{sec:evaluation}.
Theorem~\ref{hp} shows that the family of genus one normal curves
with the same theta group as $C_\phi$ is defined by taking linear
combinations of $\phi$ and $H(\phi)$.
We therefore solve for all genus one models $\phi' = \lambda \phi +
\mu H(\phi)$ with $(\lambda:\mu) \in \PP^1(K)$ and $\Jac(C_{\phi'}) \isom E$. 
We know by Theorem~\ref{xn-family}(i) that there is always at least
one such model. As explained at the end of Section~\ref{Theta groups} 
we can use the Hesse
polynomials to solve for $(\lambda:\mu)$. Notice that if $\phi$ is 
chosen to have the same invariants as a fixed Weierstrass equation
for $F$ then we can use the same $(\lambda:\mu)$ for all
$P \in F(K)$. Finally we check to see whether $C = C_{\phi'}$ is
everywhere locally soluble. If so then it represents an element
of $\Sha(E/K)[n]$. We refer to \cite{Mazur-appendix}  for some 
theoretical explanation as to why the torsors computed using 
visibility often turn out to be everywhere locally soluble.

The definition of visibility requires that the middle isomorphism
in~(\ref{eq:composite}) is the natural one
induced by the isomorphism $E[n] \isom F[n]$. However in the 
above construction we have identified $H^1(K,E[n])$ and $H^1(K,F[n])$
as parametrising the theta groups for $E[n] \isom F[n]$, first as
twists of $\Theta_E$ and then as twists of $\Theta_F$. So we are
only using the natural isomorphism if $\Theta_E \isom \Theta_F$. 
Lemma~\ref{lem:n_odd} shows that this is always true for $n$ odd,
but it can fail if $n$ is even. Nonetheless our construction defines
an explicit map $F(K)/nF(K) \to H^1(K,E)[n]$, and even though 
it may differ by a shift from the one specified by Mazur, 
it seems in numerical examples to be just as good for 
constructing elements of $\Sha(E/K)[n]$. 

The method extends to pairs of elliptic curves $E$ and $F$ 
that are reverse $n$-congruent. Instead of taking linear combinations 
of $\phi$ and its Hessian, we take linear combinations of $P(\phi)$ 
and $Q(\phi)$
where $P$ and $Q$ are the contravariants. The role of the 
Hesse polynomials is then taken by the dual Hesse polynomials.
If $\Theta$ is a theta group then we write $\Theta^\vee$ for the dual
theta group obtained by replacing the map
$\alpha : \Gm \to \Theta$ by $\lambda \mapsto \alpha(\lambda)^{-1}$. 
The analogue of Lemma~\ref{lem:n_odd} is

\begin{Lemma}
\label{lem:n_odd_dual}
If $n$ is odd and $E$ and $F$ are reverse $n$-congruent then 
$\Theta_E \isom \Theta_F^\vee$.
\end{Lemma}
\begin{Proof}
The proof of Lemma~\ref{lem:n_odd} carries over immediately.
\end{Proof}

\begin{Remark}
Bruin and Dahmen \cite{BD} have used our construction for
reverse $3$-congruent elliptic curves to show that every element of 
$\Sha(E/K)[3]$ is visible in the Jacobian of a genus $2$ curve. 
\end{Remark}

\begin{Remark}
\label{other5}
Since $-1$ is a square mod $5$ there is no difference between 
direct and reverse $5$-congruence. The case where
$E[5]$ and $F[5]$ are isomorphic as Galois modules, but the isomorphism
cannot be chosen to respect the Weil pairing, will be discussed 
further in \cite{enqI}. Our solution in this case relies on algorithms
for evaluating the covariants mentioned in Remark~\ref{othercov}.
\end{Remark}

\section{Examples}
\label{sec:examples}

We illustrate the theory in Section~\ref{sec:vis} by computing
equations for some visible elements of $\Sha(E/\Q)$.
In each case we start with a pair of $n$-congruent elliptic 
curves $E/\Q$ and $F/\Q$ taken from the paper of Cremona and 
Mazur \cite{CM}. The Mordell Weil group $F(\Q)$ is then used to 
``explain'' certain elements in the Tate-Shafarevich group $\Sha(E/\Q)$.
As the method does not rely on the existence of rational
isogenies we have deliberately chosen examples where $E/\Q$
is the only elliptic curve in its isogeny class.

Assuming a suitable elliptic curve $F/\Q$ is known, computing equations for 
elements of $\Sha(E/\Q)[n]$ using visibility is very
much faster than doing an $n$-descent. This is especially true in the
case $n=5$ where the problem of computing class groups and units often
makes $5$-descent impractical.

In presenting these examples we use the algorithms for minimising
and reducing genus one models of degrees $n=2,3,4$ described 
in \cite{minred234}. In other words we make changes of co-ordinates
so that the models considered have small integer coefficients.
We do not record the changes of 
co-ordinates used, as these may be recovered using the algorithms in 
\cite{bq}, \cite{testeqtc}, \cite{improve4}, as implemented in 
our {\sf MAGMA} function {\tt IsEquivalent}. In 
the examples with $n > 2$ the models $\phi_i$ are the ones returned by 
our {\sf MAGMA} function {\tt GenusOneModel(n,P)}.


\subsection{An example of $\Sha(E/\Q)[2]$} 
There is no distinction between direct and reverse $2$-congruence.
The formulae in Section \ref{famec} show that the elliptic curves
\begin{align*} E = {\rm 571a1} : \quad 
& y^2 + y = x^3 - x^2 - 929 x - 10595 \\
F  = {\rm 571b1}: \quad
& y^2 + y = x^3 + x^2 - 4 x + 2
\end{align*}
are $2$-congruent. 
We have $E(\Q)=0$ and $F(\Q) \isom \Z^2$ 
generated by $P_1 = (0,1)$ and $P_2 = (1 , 0)$. 
Mapping $F \to \PP^1$ via the complete
linear system $|0_F + P|$ for $P = 0_F, P_1, P_2, P_1+P_2$
we obtain binary quartics
\begin{align*}
\phi_1 &=  4 x^3 z + 16 x^2 z^2 + 4 x z^3 + z^4 \\
\phi_2 &=  x^4 + 4 x^3 z + 4 x^2 z^2 - 12 x z^3 + 4 z^4 \\
\phi_3 &=  x^4 + 4 x^3 z - 2 x^2 z^2 - 8 x z^3 + 9 z^4 \\
\phi_4 &=  x^4 - 8 x^3 z + 10 x^2 z^2 + 4 x z^3 + z^4
\end{align*}
with invariants $c_4 = 3328$, $c_6 = -202240$
and $\Delta = -2^{12} \cdot 571$.
Let $\D,\cc_4,\cc_6$ be the Hesse polynomials with coefficients 
evaluated at $c_4, c_6$.
The binary form \[ \cc_4(\lambda,\mu)^3 - j(E) \Delta \D(\lambda,\mu)^2 =0\]
has a unique $\Q$-rational root at $(\lambda:\mu) = (-116:1)$. 
Solving for $d \in \Q^\times$ with  
$\cc_4(-166,1) = d^2 c_4(E)$ and $\cc_6(-166,1) = d^3 c_6(E)$ we find
$d \in 3 (\Q^\times)^2$.
We therefore compute $3(-116 \phi_i + H(\phi_i))$ for $i=1,2,3,4$ and then 
minimise and reduce to obtain
\begin{align*}
\psi_1 &=  -4 x^4 - 60 x^3 z - 232 x^2 z^2 - 52 x z^3 - 3 z^4 \\
\psi_2 &=  -11 x^4 - 68 x^3 z - 52 x^2 z^2 + 164 x z^3 - 64 z^4 \\
\psi_3 &=  -15 x^4 - 52 x^3 z + 38 x^2 z^2 + 144 x z^3 - 115 z^4 \\
\psi_4 &=  -19 x^4 + 112 x^3 z - 142 x^2 z^2 - 68 x z^3 - 7 z^4.
\end{align*}
Each of these binary quartics $\psi_i$ 
has discriminant $-2^{12} \cdot 571$, and defines a curve
that is local solubility at $p= 2,571$.
Real solubility is automatic since the discriminant
is negative. The binary quartic $\psi_3$ has a rational root
at $(x:z) = (1:1)$. It follows that the $\psi_i$ define a subgroup of
$\Sha(E/\Q)[2]$ isomorphic to $(\Z/2\Z)^2$ with $\psi_3$ corresponding
to the identity. The fact $\phi_1$ represents the identity in 
$H^1(\Q,F[2])$ whereas $\psi_3$ represents the identity in $H^1(\Q,E[2])$
shows that in this example $\Theta_E$ and $\Theta_F$ are not isomorphic.

\subsection{Examples of $\Sha(E/\Q)[3]$} 
We give two examples, one where the elliptic curves $E$ and $F$
are directly $3$-congruent, and another where they are reverse
$3$-congruent.
The formulae in Section \ref{famec} show that the elliptic curves
\begin{align*} E = {\rm 2006e1} : \quad 
& y^2 + x y = x^3 + x^2 - 58293654 x - 171333232940 \\
F  = {\rm 2006d1}: \quad
& y^2 + x y = x^3 + x^2 - 88 x + 284 
\end{align*}
are directly $3$-congruent. 
We have $E(\Q)=0$ and $F(\Q) \isom \Z^2$ 
generated by $P_1 = (-10,22)$ and $P_2 = (2 , 10)$. 
Embedding $F \subset \PP^2$ via the complete
linear system $|2.0_F + P|$ for $P = P_1, P_2, P_1 + P_2, P_1 - P_2$
we obtain ternary cubics
\small
\begin{align*}
\phi_1 &=  x^2 y - 2 x^2 z + x y^2 - x y z - x z^2 - 2 y^3 + y^2 z 
       + 5 y z^2 + 2 z^3  \\
\phi_2 &= -x^2 y - x y^2 - 5 x y z + x z^2 + 2 y^2 z + 9 y z^2 - z^3  \\
\phi_3 &= -x^2 y + 2 x y^2 - 7 x y z + x z^2 - y^2 z + 6 y z^2 - z^3  \\
\phi_4 &= x^3 + 3 x^2 y + 2 x^2 z + x y^2 + x y z - 2 x z^2 - y^3 
      + 2 y^2 z + y z^2 - 2 z^3 
\end{align*}
\normalsize
with the same invariants
$c_4 = 4249$, $c_6 = -277181$ and $\Delta = -2^2 \cdot 17^2 \cdot 59$ 
as $F$. 
Let $\D,\cc_4,\cc_6$ be the Hesse polynomials with coefficients evaluated
at $c_4, c_6$. The binary form 
\[\cc_4(\lambda,\mu)^3 - j(E) \Delta \D(\lambda,\mu)^3 = 0 \] 
has a unique $\Q$-rational 
root at $(\lambda:\mu) = (521:9)$. We therefore compute
$521 \phi_i + 9 H(\phi_i)$ for $i=1,2,3,4$ and then minimise and reduce
to obtain 
\small
\begin{align*}
\psi_1 & = 9 x^3 - 16 x^2 y + 5 x^2 z + 38 x y^2 + 129 x y z 
       + 6 x z^2 + 59 y^3 - 81 y^2 z - 58 y z^2 - 124 z^3 \\
\psi_2 & =     9 x^3 + 43 x^2 y - 27 x^2 z + 75 x y^2 
       + 53 x y z + 92 x z^2 - 4 y^3 + 
        75 y^2 z + 2 y z^2 + 124 z^3 
\end{align*}
\begin{align*}
\psi_3 & =     9 x^3 + 43 x^2 y - 27 x^2 z + 27 x y^2 
       + 85 x y z - 43 x z^2 + 74 y^3 + 
        74 y^2 z - 58 y z^2 - 92 z^3 \\
\psi_4 & =     43 x^3 + 38 x^2 y + 22 x^2 z - 48 x y^2 
       - 43 x y z + 65 x z^2 + 11 y^3 - 
        5 y^2 z + 113 y z^2 + 50 z^3 .
\end{align*}
\normalsize
Each of these ternary cubics has the same invariants as $E$, and 
defines a curve that is locally soluble at $p=2,17,59$ 
(the bad primes of $E$).
It follows that the $\psi_i$ define the inverse pairs of non-zero elements 
in a subgroup of $\Sha(E/\Q)[3]$ isomorphic to $(\Z/3\Z)^2$.

Our second example is similar.
The formulae in Section \ref{famec} show that the elliptic curves
\begin{align*} E = {\rm 2541d1 } : \quad 
& y^2 + y = x^3 - x^2 - 180572 x - 26845765 \\
F  = {\rm 2541c1}: \quad
& y^2 + x y + y = x^3 + x^2 + 3 x + 12 
\end{align*}
are reverse $3$-congruent.
We have $E(\Q)=0 $ and $F(\Q) \isom \Z^2$ generated by
$P_1 = (-2 , 2)$ and $P_2 = (0, 3)$.
Embedding $F \subset \PP^2$ via the complete
linear system $|2.0_F + P|$ for $P =  P_1, P_2, P_1 + P_2, P_1 - P_2$
we obtain ternary cubics
\small
\begin{align*}
\phi_1 &= -x^2 z + x y^2 - x y z + x z^2 + 2 y^2 z + y z^2 - 6 z^3 \\
\phi_2 &= -x^2 z + x y^2 + x y z + x z^2 - y^2 z + 6 y z^2 \\
\phi_3 &= -x^2 y + x y^2 + x y z + 2 x z^2 + 2 y^2 z - 3 y z^2 + z^3 \\
\phi_4 &= -x^2 y + x y z + x z^2 + y^3 + 2 y^2 z - 2 y z^2 + 2 z^3
\end{align*}
\normalsize
with the same invariants $c_4 = -143$, $c_6 = -9449$ and 
$\Delta = -3^2 \cdot 7^2 \cdot 11^2$ as $F$. 
Let $\DD,\cdual_4,\cdual_6$ be the dual Hesse polynomials with coefficients 
evaluated at $c_4, c_6$. The binary form
 \[\cdual_4(\xi,\eta)^3 
-  1728 j(E) \Delta^2 \DD(\xi,\eta)^3 = 0\] has a unique $\Q$-rational 
root at $(\xi:\eta) = (-55:1)$. We therefore compute
$-55 P(\phi_i) + Q(\phi_i)$ for $i=1,2,3,4$ and then minimise and reduce
to obtain
\small
\begin{align*}
\psi_1 & =  -x^3 - 3 x^2 y - 7 x^2 z - 14 x y^2 + 8 x y z 
    + 13 x z^2 - y^3 + 26 y^2 z + 2 y z^2 + 70 z^3 \\
\psi_2 & = -3 x^3 - 14 x^2 y - 5 x^2 z - x y^2 + 4 x y z - 15 x z^2 
    - 5 y^3 + 30 y^2 z - 16 y z^2 - 26 z^3 \\
\psi_3 & = 3 x^3 + 7 x^2 y - 4 x^2 z + 3 x y^2 + 28 x y z + 25 x z^2 
    - 9 y^3 - 5 y^2 z + 6 y z^2 + 35 z^3 \\
\psi_4 & = x^3 + 7 x^2 y - 12 x^2 z + 9 x y^2 + 10 x y z + 37 x z^2 
    - 4 y^3 + 8 y^2 z + 2 y z^2 + 35 z^3.
\end{align*}
\normalsize
Each of these ternary cubics has the same invariants as $E$, and defines a 
curve that is locally soluble at $p=3,7,11$ (the bad primes of $E$).
It follows that the $\psi_i$ define the inverse pairs of non-zero 
elements in a subgroup of $\Sha(E/\Q)[3]$ isomorphic to 
$(\Z/3\Z)^2$.

\subsection{Examples of $\Sha(E/\Q)[4]$}
We give two examples, one where the elliptic curves $E$ and $F$
are directly $4$-congruent, and another where they are reverse
$4$-congruent. 
The formulae in Section \ref{famec} show that the elliptic curves
\begin{align*} E = {\rm 2045b1} : \quad 
& y^2 + x y = x^3 - x^2 - 5470 x - 862675 \\
F = {\rm 4090b1}: \quad
& y^2 + x y = x^3 + x^2 + 7 x + 37 
\end{align*}
are directly $4$-congruent. 
We have $E(\Q)=0$ and $F(\Q) \isom \Z^2$ 
generated by $P_1 = (2,7)$ and $P_2 = (18 , 71)$. 
Embedding $F \subset \PP^3$ via the complete
linear system $|3.0_F + P|$ for 
$P = P_1,P_2,P_1+P_2,P_1+2 P_2,P_1-P_2,2 P_1+P_2$
we obtain quadric intersections
\small
\begin{align*}
\phi_1 &= \begin{pmatrix}
        x_1 x_4 - x_2 x_3 - x_2 x_4 + x_3^2 - x_3 x_4 + 2 x_4^2 \\
        x_1 x_3 + x_1 x_4 + x_2^2 - x_2 x_3 + x_3^2 - 7 x_3 x_4 - 4 x_4^2
    \end{pmatrix} \\
\phi_2 &= \begin{pmatrix}
        x_1 x_3 + x_2^2 + x_2 x_4 - x_3^2 - 2 x_3 x_4 - 2 x_4^2 \\
        x_1 x_3 + x_1 x_4 + x_2^2 - x_2 x_3 + 3 x_3^2 - x_3 x_4 - 2 x_4^2
    \end{pmatrix} \\
\phi_3 &= \begin{pmatrix}
        x_1 x_4 - x_2 x_3 + x_2 x_4 + 3 x_4^2 \\
        x_1 x_2 + x_1 x_4 - 8 x_2 x_4 + x_3^2 + 4 x_4^2
    \end{pmatrix} \\
\phi_4 &= \begin{pmatrix} 
        x_1 x_3 - x_2 x_4 + x_3^2 - x_3 x_4 + x_4^2 \\
        x_1 x_2 - x_1 x_3 - 2 x_2 x_3 + x_2 x_4 + 3 x_4^2
    \end{pmatrix} \\
\phi_5 &= \begin{pmatrix}
        x_1 x_2 + x_1 x_4 - 2 x_2 x_3 + 2 x_2 x_4 + x_3^2 - 2 x_4^2 \\
        -x_1 x_4 + 2 x_2^2 + x_2 x_3 + 3 x_2 x_4 + x_4^2
    \end{pmatrix} \\
\phi_6 &= \begin{pmatrix}
        x_1 x_3 + x_2 x_3 + 3 x_2 x_4 + x_3^2 + x_4^2 \\
        x_1 x_4 + x_2^2 - x_2 x_3 - 3 x_3 x_4 - x_4^2
    \end{pmatrix} 
\end{align*}
\normalsize
with the same invariants 
$c_4 = -311$, $c_6 = -29573$ 
and $\Delta = -2^8 \cdot 5 \cdot 409$ as $F$.
Let $\D,\cc_4,\cc_6$ be the Hesse polynomials with coefficients evaluated
at $c_4, c_6$. The binary form 
 \[\cc_4(\lambda,\mu)^3 
- j(E) \Delta \D(\lambda,\mu)^4 = 0\] has a unique $\Q$-rational 
root at $(\lambda:\mu) = (5:1)$. 
We therefore compute $\phi'_i = 5 \phi_i + H(\phi_i)$ for $i=1, \ldots ,6$ 
and then minimise and reduce to obtain
\small
\begin{align*}
\psi_1 &= \begin{pmatrix}
        x_1 x_2 + 2 x_1 x_4 - x_2 x_3 - 4 x_2 x_4 + x_3^2 + x_3 x_4 + x_4^2 \\
        x_1^2 + 2 x_1 x_2 + x_1 x_3 + 3 x_1 x_4 + 7 x_2^2 - x_2 x_3 + 2 x_3^2 - 4 x_3 x_4 - 
            2 x_4^2
    \end{pmatrix} 
\end{align*}
\begin{align*}
\psi_2 &= \begin{pmatrix}
        2 x_1 x_2 + x_1 x_3 + x_1 x_4 + x_2 x_3 - 2 x_2 x_4 + 3 x_3^2 + 2 x_3 x_4 + 4 x_4^2 \\
        x_1^2 - x_1 x_2 - x_2^2 - 5 x_2 x_3 + 4 x_2 x_4 - 2 x_3^2 + x_3 x_4
    \end{pmatrix} \\
\psi_3 &= \begin{pmatrix}
        x_1 x_3 + 3 x_1 x_4 + x_2^2 - x_2 x_3 + x_3^2 - 3 x_3 x_4 + x_4^2 \\
        x_1^2 + 4 x_1 x_2 + 2 x_1 x_3 - 6 x_1 x_4 - x_2 x_3 - 2 x_2 x_4 - 2 x_3^2 - 2 x_3 x_4 
            + 2 x_4^2
    \end{pmatrix} \\ 
\psi_4 &= \begin{pmatrix}
        2 x_1 x_2 + 2 x_1 x_3 + x_1 x_4 + x_2^2 + x_2 x_3 + x_3^2 + x_3 x_4 + 2 x_4^2 \\
        x_1^2 + x_1 x_2 - 3 x_1 x_3 - 4 x_1 x_4 + 2 x_2^2 + 3 x_2 x_4 + 3 x_3^2 - 2 x_3 x_4 -
            x_4^2
    \end{pmatrix} \\
\psi_5 &= \begin{pmatrix}
        x_1^2 + x_1 x_2 + x_1 x_3 + x_1 x_4 - x_2^2 - 2 x_2 x_3 + x_3^2 + x_3 x_4 + 2 x_4^2 \\
        x_1^2 - 4 x_1 x_2 - 2 x_2^2 - 4 x_2 x_3 + 5 x_2 x_4 - 3 x_3^2 + 2 x_3 x_4 - x_4^2
    \end{pmatrix} \\
\psi_6 &= \begin{pmatrix}
        x_1 x_2 + 3 x_1 x_3 + x_2^2 + x_2 x_4 + x_3^2 - 2 x_3 x_4 + 5 x_4^2 \\
        x_1^2 - x_1 x_2 - x_1 x_3 + 7 x_1 x_4 + x_2 x_3 + 3 x_2 x_4 + x_3^2 + 3 x_3 x_4 - x_4^2
    \end{pmatrix} . 
\end{align*}
\normalsize
Each of these quadric intersections has the same invariants as $E$, and
defines a curve that is locally soluble at $p=5,409$ (the bad primes of $E$).
Solubility over the reals is automatic since the discriminant is negative. 
Repeating for $P = 0_F$ we obtain a quadric intersection
equivalent to the one defining $E \subset \PP^3$ embedded by $|4.0_E|$. 
So in this example $\Theta_E$ and $\Theta_F$ are 
isomorphic. It follows that the $\psi_i$ 
define the inverse pairs 
of elements of order~$4$ in a subgroup of $\Sha(E/\Q)[4]$ isomorphic to 
$(\Z/4\Z)^2$.

Our second example is similar.
The formulae in Section \ref{famec} show that the elliptic curves
\begin{align*} E = {\rm 1309a1 } : \quad 
& y^2 + y = x^3 - 406957 x - 99924251 \\
F  = {\rm 1309b1}: \quad
& y^2 + y = x^3 - x^2 - 22 x + 52 
\end{align*}
are reverse $4$-congruent.
We have $E(\Q)=0 $ and $F(\Q) \isom \Z^2$ generated by
$P_1 = (16 , 59)$ and $P_2 = (-1, 8)$.
Embedding $F \subset \PP^3$ via the complete
linear system $|3.0_F + P|$ for 
$P = P_1,P_2,P_1+P_2,P_1+2 P_2,P_1-P_2,2 P_1+P_2$
we obtain quadric intersections
\small
\begin{align*}
\phi_1 &= \begin{pmatrix}
        x_1 x_3 + x_1 x_4 + x_2 x_4 - 2 x_3 x_4 + x_4^2 \\
        x_1 x_4 + x_2^2 + x_2 x_3 - x_2 x_4 - 2 x_3^2
    \end{pmatrix} \\
\phi_2 &= \begin{pmatrix}
        x_1 x_3 + x_2 x_3 + x_2 x_4 + 2 x_3 x_4 \\
        x_1 x_4 + x_2^2 - 3 x_2 x_4 + x_3^2 + x_3 x_4 - 2 x_4^2
    \end{pmatrix} 
\end{align*}
\begin{align*}
\phi_3 &= \begin{pmatrix}
        x_1 x_3 + x_1 x_4 - x_2 x_4 + x_3^2 - x_3 x_4 - 2 x_4^2 \\
        x_1 x_3 + x_2^2 - x_2 x_4 + 3 x_3 x_4 - 2 x_4^2
    \end{pmatrix}  \\
\phi_4 &= \begin{pmatrix}
        x_1 x_2 + x_1 x_3 + x_2^2 + x_2 x_4 - x_3^2 - x_4^2 \\
        x_1 x_2 + x_2 x_3 + x_2 x_4 + 3 x_3 x_4 + x_4^2
    \end{pmatrix} \\ 
\phi_5 &= \begin{pmatrix}
        x_1 x_4 + x_2 x_3 + x_2 x_4 - x_3 x_4 + x_4^2 \\
        x_1 x_2 + 3 x_2 x_3 - 2 x_2 x_4 + x_3^2 + 3 x_3 x_4 + 2 x_4^2
    \end{pmatrix} \\
\phi_6 &= \begin{pmatrix}
        x_1 x_3 + x_2^2 - x_3^2 - x_3 x_4 - x_4^2 \\
        x_1 x_2 + x_1 x_3 + 2 x_2 x_3 - x_2 x_4 + 2 x_3 x_4 + x_4^2
    \end{pmatrix} 
\end{align*}
\normalsize
with the same invariants 
$c_4 =  1072$, $c_6 = -38744$
and $\Delta = -7^2 \cdot 11 \cdot 17^2$ as $F$.
Let $\DD,\cdual_4,\cdual_6$ be the dual Hesse polynomials with 
coefficients evaluated at $c_4, c_6$. The binary form 
 \[\cdual_4(\xi,\eta)^3 
- 1728^2 j(E) \Delta^3 \DD(\xi,\eta)^4 = 0\] has a unique $\Q$-rational 
root at $(\xi:\eta) = (35:1)$. 
We therefore compute $35 P(\phi_i) + Q(\phi_i)$ for $i=1, \ldots, 6$
and then minimise and reduce to obtain
\small
\begin{align*}
\psi_1 &= \begin{pmatrix}
        x_1^2 + 2 x_1 x_2 + 4 x_1 x_3 + x_1 x_4 + 2 x_2^2 + 7 x_2 x_3 + x_2 x_4 + 2 x_3^2 - 
            8 x_3 x_4 + 7 x_4^2 \\
        2 x_1 x_2 + x_1 x_3 + x_1 x_4 + x_2^2 + 2 x_2 x_3 + 13 x_3^2 - 2 x_3 x_4 + 4 x_4^2
    \end{pmatrix} \\
\psi_2 &= \begin{pmatrix}
        x_1 x_3 + x_1 x_4 + x_2^2 - 4 x_2 x_3 - 4 x_3^2 - 17 x_3 x_4 - 8 x_4^2 \\
        x_1^2 + x_1 x_4 + x_2 x_3 - 3 x_2 x_4 + x_3^2 - 4 x_3 x_4 + 20 x_4^2
    \end{pmatrix} \\
\psi_3 &= \begin{pmatrix}
        x_1^2 + x_1 x_3 + x_2^2 + x_2 x_3 - x_2 x_4 - x_3^2 - 4 x_3 x_4 + 3 x_4^2 \\
        5 x_1 x_2 + 3 x_1 x_3 + 3 x_1 x_4 + 2 x_2^2 + 2 x_2 x_3 + 4 x_2 x_4 - 7 x_3^2 - 
            4 x_3 x_4 - 8 x_4^2
    \end{pmatrix} \\
\psi_4 &= \begin{pmatrix}
        x_1^2 + x_1 x_2 + 2 x_1 x_3 + 5 x_1 x_4 + x_2^2 + 3 x_2 x_3 + 6 x_2 x_4 + 2 x_3^2 - 
            2 x_3 x_4 - 7 x_4^2 \\
        2 x_1^2 - 2 x_1 x_3 + 6 x_1 x_4 + 2 x_2^2 + x_2 x_3 + 7 x_2 x_4 + 2 x_3^2 - 5 x_3 x_4
            + 4 x_4^2
    \end{pmatrix} \\
\psi_5 &= \begin{pmatrix}
        4 x_1 x_2 + 4 x_1 x_3 + x_1 x_4 - 6 x_2 x_3 - 4 x_2 x_4 + x_3^2 - 3 x_3 x_4 + x_4^2 \\
        x_1^2 + x_1 x_2 - x_1 x_3 + 2 x_1 x_4 + 7 x_2^2 - 5 x_2 x_3 - 4 x_2 x_4 + x_3^2 + 
            x_3 x_4 + 2 x_4^2
    \end{pmatrix} \\
\psi_6 &= \begin{pmatrix}
        3 x_1 x_3 + 6 x_1 x_4 + x_2^2 + x_3^2 - x_3 x_4 + 9 x_4^2 \\
        x_1^2 + 3 x_1 x_2 - 6 x_1 x_3 - 10 x_1 x_4 + 2 x_2 x_3 + 3 x_2 x_4 - x_3^2 - 5 x_3 x_4
            + 2 x_4^2
    \end{pmatrix} .
\end{align*}
\normalsize
Each of these quadric intersections has the same invariants as $E$, and 
defines a curve that is locally soluble at $p=7,11,17$ (the bad primes of $E$).
Solubility over the reals is automatic since the discriminant is negative. 
Repeating for $P = 2P_2 $ we obtain a quadric intersection
equivalent to the one defining $E \subset \PP^3$ embedded by $|4.0_E|$. 
So in this example the theta groups $\Theta_E$ and $\Theta_F^\vee$ differ 
by a $2$-torsion element in $H^1(\Q,E[4]) \isom H^1(\Q,F[4])$. 
It follows that the $\psi_i$ define the inverse pairs 
of elements of order $4$ in a subgroup of $\Sha(E/\Q)[4]$ isomorphic to 
$(\Z/4\Z)^2$.

\subsection{An example of $\Sha(E/\Q)[5]$} 
The formulae in Section \ref{famec} show that the elliptic curves
\begin{align*} E = {\rm 1058d1} : \quad 
& y^2 + x y = x^3 - x^2 - 332311 x - 73733731 \\
F = {\rm 1058c1}: \quad
& y^2 + x y + y = x^3 + 2
\end{align*}
are directly $5$-congruent. 
We have $E(\Q)=0$ and $F(\Q) \isom \Z^2$ 
generated by $P_1 = (-1,1)$ and $P_2 = (0 , 1)$. 
Embedding $F \subset \PP^4$ via the complete
linear system $|4.0_F + P|$ for $P= P_1$, i.e.
\[ (x,y) \mapsto  (x^2+1 : -y-1 : x : (y+1)/(x+1) : 1 ), \]
we obtain (using the algorithm in \cite{g1pf}) a genus one model
\small
\[ \phi_1 = \begin{pmatrix} 
    0  & -x_1 + x_3 - x_5   &    x_4   &   x_2 + x_4   &    -x_4  \\
   &   0   &    x_2 + x_5   &   -x_1 + x_5  &    -x_3  \\
             &      &     0    &    x_3   &     x_5  \\
      &  -  &     &    0    &    0  \\
             &      &      &         &   0
\end{pmatrix} \]
\normalsize
with the same invariants $c_4= -23$, $c_6 = -1909$ 
and $\Delta = - 2^2 \cdot 23^2$ as $F$.
Let $\D,\cc_4,\cc_6$ be the Hesse polynomials with coefficients evaluated
at $c_4, c_6$. The binary form 
 \[\cc_4(\lambda,\mu)^3 
- j(E) \Delta \D(\lambda,\mu)^5 = 0 \] has a unique $\Q$-rational 
root at $(\lambda:\mu) = (-23:1)$. 
We compute the Hessian $H(\phi_1)$ using the algorithm in 
Section~\ref{sec:evaluation} and find it has entries 
\small
\begin{align*}
 H_{1 2} &= x_1 - 61 x_3 - 35 x_5 & H_{2 4} &= x_1 - 12 x_2 + 12 x_3 + 47 x_5 \\
 H_{1 3} &= 12 x_1 - 12 x_2 + 36 x_3 - 13 x_4 - 60 x_5
&  H_{2 5} &= -12 x_1 - 12 x_2 + 25 x_3 - 24 x_4 - 36 x_5 \\ 
 H_{1 4} &= -x_2 - 12 x_3 - 37 x_4 - 12 x_5 & 
 H_{3 4} &= -24 x_2 + 35 x_3 - 24 x_4 - 48 x_5  \\
 H_{1 5} &= 12 x_2 - 12 x_3 - 11 x_4 + 12 x_5 & H_{3 5} &= -12 x_3 - x_5 \\
 H_{2 3} &= 12 x_1 + 23 x_2 - 12 x_3 + 72 x_4 + 47 x_5 & 
H_{4 5} &= -24 x_2 + 12 x_3 - 12 x_5 
\end{align*}
\normalsize
Then we compute $\phi'_1 = -23 \phi_1 + H(\phi_1)$ and make a 
transformation of
the form $[A,\lambda I_5] : \phi'_1 \mapsto \lambda A \phi'_1 A^T$ for 
some $\lambda \in \Q^\times$ and $A \in \GL_5(\Q)$
to obtain
\small
\[ \psi_1 = \begin{pmatrix} 
0 &  -x_2 &  x_5  & -x_2 + x_3 + x_5 &  -x_4 - 2 x_5 \\
 &  0 &  2 x_4  & -x_1 - 2 x_5 &  -x_3 + 2 x_4 - x_5 \\
 &  &  0 &  -2 x_3 &  x_1 + x_3 + x_4 \\
 & - &  &  0 &  
-x_2 + 4 x_3 + 4 x_4 - 2 x_5 \\
 &  &  &  &  0
\end{pmatrix}. \]
\normalsize
This genus one model has the same invariants as $E$, and its $4$ by $4$ 
Pfaffians define a curve $C \subset \PP^4$ that is locally soluble at
$p=2,23$ (the bad primes of $E$). It follows that $C$ represents
a non-trivial element of $\Sha(E/\Q)[5]$. We note that, unlike the 
examples in \cite{jems}, the elliptic curve $E$ does not admit 
any $\Q$-rational isogenies of degree $5$.
Repeating for $P = r_1 P_1 + r_2 P_2$ for $0 \le r_1,r_2 \le 4$ 
we find equations
for all elements in a subgroup of $\Sha(E/\Q)[5]$ isomorphic to 
$(\Z/5\Z)^2$. The full list of equations, together with similar
examples for other elliptic curves $E/\Q$ of small conductor, 
may be found on the author's website.



\end{document}